\RequirePackage{fix-cm}

\documentclass[USenglish]{article}

\usepackage[big,online]{dgruyter}	
\usepackage{helvet}

\usepackage{sfmath}
\usepackage[expansion=false]{microtype} 
\usepackage[numbers,square,sort&compress]{natbib}

\makeatletter
\AtBeginDocument{\sbox\dg@wordmark{\includegraphics[height=2.1mm]{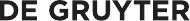}}}
\makeatother

\theoremstyle{dgthm}
\newtheorem{theorem}{Theorem}[section]
\newtheorem{corollary}[theorem]{Corollary}
\newtheorem{proposition}[theorem]{Proposition}
\newtheorem{lemma}[theorem]{Lemma}
\newtheorem{assump}{Assumption}

\theoremstyle{dgdef}

\newtheorem{remark}{Remark}

    \usepackage{mathtools, amssymb, amsfonts, cancel, array, mathrsfs, dsfont}
    \usepackage{xcolor, tikz}
    \usetikzlibrary{arrows.meta}
    \usepackage{hyperref}
    \usepackage{multirow}

    \DeclareFontFamily{U}{rsfs}{\skewchar\font127 }
    \DeclareFontShape{U}{rsfs}{m}{n}{
       <-6.5>  rsfs5
       <6.5-8> rsfs7
       <8->    rsfs10
    }{}

    \let\pa\partial
    
    \newcommand{\cA}{\mathcal{A}}

    \newcommand{\cE}{\mathcal{E}}

    \newcommand{\cL}{\mathcal{L}}
    \newcommand{\cP}{\mathcal{P}}
    \newcommand{\cR}{\mathcal{R}}
    \newcommand{\cT}{\mathcal{T}}
    \newcommand{\cV}{\mathcal{V}}
    
    \newcommand{\sD}{\mathscr{D}}
    
    \newcommand{\ba}{\mathbf{a}}
    \newcommand{\be}{\mathbf{e}}
    
    \newcommand{\bk}{\mathbf{k}}
    \newcommand{\bn}{\mathbf{n}}
    \newcommand{\bp}{\mathbf{p}}
    
    \newcommand{\bs}{\mathbf{s}}
    \newcommand{\bu}{\mathbf{u}}
    \newcommand{\bv}{\mathbf{v}}
    \newcommand{\bw}{\mathbf{w}}
    \newcommand{\bx}{\mathbf{x}}
    \newcommand{\by}{\mathbf{y}}
    \newcommand{\bz}{\mathbf{z}}
    \newcommand{\bzero}{\mathbf{0}}
    
    \newcommand{\bphi}{\boldsymbol{\phi}}

    \newcommand{\bC}{\mathbf{C}}
    \newcommand{\bF}{\mathbf{F}}
    \newcommand{\bH}{\mathbf{H}}
    \newcommand{\bI}{\mathbf{I}}
    
    \newcommand{\bL}{\mathbf{L}}
    \newcommand{\bM}{\mathbf{M}}
    \newcommand{\bP}{\mathbf{P}}
    
    \newcommand{\bS}{\mathbf{S}}
    \newcommand{\bV}{\mathbf{V}}
    \newcommand{\bW}{\mathbf{W}}

    \newcommand{\dL}{\mathds{L}}

    \newcommand{\rJ}{\mathrm{J}}

    \newcommand{\fC}{\mathfrak{C}}

    \newcommand{\uD}{\underline{D}}
    
    \newcommand{\bbR}{\mathbb{R}}
    \newcommand{\bbZ}{\mathbb{Z}}

    \newcommand{\wt}{\widetilde}
    \newcommand{\wh}{\widehat}

    \newcommand{\ophi}{\overline{\phi}}
    \newcommand{\oc}{e^{c}}
	\newcommand{\omu}{e^{\mu}}
	\newcommand{\obu}{\be^{u}}
	\newcommand{\op}{e^{p}}

    \newcommand{\wbu}{\widetilde{\mathbf{u}}}
    \renewcommand{\wp}{\widetilde{p}}
    \newcommand{\wc}{\widetilde{c}}
    \newcommand{\wmu}{\widetilde{\mu}}
    
    \renewcommand{\div}{\operatorname{div}}
    \newcommand{\divG}{\mathrm{div}_\G}
    \newcommand{\DeltaG}{\Delta_\G}
    \renewcommand{\deg}{\operatorname{deg}}
    
    \newcommand{\G}{\Gamma}
    \newcommand{\dst}{\operatorname{d}}
    \newcommand{\nablaG}{\nabla_\Gamma}
    
    \newcommand{\dt}{\,\mathrm{d}t}
    \newcommand{\ds}{\,\mathrm{d}s}
    \newcommand{\dA}{\,\mathrm{d}A}
    \newcommand{\dx}{\,\mathrm{d}x}
    \newcommand{\dz}{\,\mathrm{d}\mathbf{z}}
    \newcommand{\ddt}{\frac{\mathrm{d}}{\mathrm{d}t}}
    \newcommand{\mean}[1]{\langle #1 \rangle_\G}
    \newcommand{\zm}[1]{#1^\#}
    \newcommand{\tr}{\operatorname{tr}}
    
    \newcommand{\rev}[1]{{\color{blue}#1}} 
    
    \newcommand{\Csob}{C_{\mathrm{Sob}}}
    
    \newcommand{\Cint}{C_{\mathrm{int}}}

\begin{document}

\title{Analysis of TraceFEM for the Surface Stokes--Cahn--Hilliard Model}
\runningtitle{Analysis of TraceFEM for Surface Stokes--Cahn--Hilliard}

\author[1]{Maxim Olshanskii}
\author*[2]{Mansur Shakipov}

\runningauthor{M. A. Olshanskii and M. Shakipov}

\affil[1]{\protect\raggedright 
University of Houston, Department of Mathematics, Houston, USA, e-mail: maolshanskiy@uh.edu}
\affil[2]{\protect\raggedright 
University of Houston, Department of Mathematics, Houston, USA, e-mail: mshakipov@uh.edu}

    
\dedication{the memory of  Prof.~Dr.~Yuri Alexeevich Kuznetsov.}
	
\abstract{
    We analyze a trace finite element method (TraceFEM) for a surface fluid with a stabilization term that filters out Killing vector fields, and use it to discretize the surface Stokes--Cahn--Hilliard problem in space only. Assuming no geometric error, we show that the discrete inf-sup condition that our method requires is equivalent to the usual condition found in the literature and known to hold for trace Taylor--Hood pairs of arbitrary order in intrinsic dimension two. Inspired by a continuous dependence estimate at the PDE level, we derive error estimates for our scheme in the energy norms, which in turn lead to estimates in $L^2$-norms and, for two-dimensional surfaces, $L^\infty$-norm for the phase variable. To that end, we introduce a $p$-dependent scaling for the normal derivative volume stabilization, and use it to prove new approximation and inverse estimates in the $L^p$- and $W^{1,p}$-based setting ($p \in [1,\infty]$ for approximation and $p \in [2,\infty]$ for inverse estimates). Finally, we corroborate the analysis with numerical experiments.
}


\maketitle

\section{Introduction}

Lateral phase separation in multicomponent lipid membranes produces domains whose evolution depends on both diffusion and tangential fluid motion. Surface Navier--Stokes--Cahn--Hilliard models have been used to study this coupling~\cite{NitschkeVoigtWensch2012,PalzhanovZhiliakovQuainiOlshanskii2021}. Simulations with coupled phase-field and flow models have also been compared with microscopy measurements of domain coarsening in lipid vesicles~\cite{WangPalzhanovQuainiOlshanskiiMajd2022lipid,TangWiseLowengrubGuo2026}.
At the scales of lipid-domain coarsening, representative material parameters yield a small surface Reynolds number, which motivates neglecting fluid inertia. Inertia-free descriptions of surface flow have also been used in models of deformable fluid membranes~\cite{ArroyoDeSimone2009,TorresSanchezMillanArroyo2019,PorrmannVoigt2024}. Here we consider the surface Stokes--Cahn--Hilliard system studied in Part~I~\cite{OlshanskiiShakipov2026}, which arises as a formal inertia-free reduction of the surface Navier--Stokes--Cahn--Hilliard equations.
In this system, capillary forces drive tangential flow, which in turn advects the phase field. Its numerical approximation must account for this nonlinear coupling, the tangentiality of the velocity, and the Killing fields in the kernel of the surface viscous operator.

In the companion paper~\cite{OlshanskiiShakipov2026},  we established well-posedness of the surface Stokes--Cahn--Hilliard system and introduced a structure-preserving trace finite element method (TraceFEM). The method uses traces of continuous bulk finite element functions on an implicitly represented surface~\cite{OlshanskiiReuskenGrande2009,OlshanskiiReusken2017}. Normal derivative volume stabilization controls the bulk degrees of freedom independently of how the surface intersects the background mesh~\cite{GrandeLehrenfeldReusken2018,BurmanHansboLarsonMassing2018}. A further penalty controls the Killing fields, which form the kernel of the surface viscous operator. Applying the same Killing projection in the transport and capillary terms preserves the cancellation underlying the energy balance. The resulting discretization conserves the total phase and satisfies a discrete energy identity.

The present paper studies stability and convergence of the spatial semidiscretization. We consider closed hypersurfaces $\G\subset\bbR^{d+1}$ of intrinsic dimension $d=2,3$, constant positive viscosity and mobility, and regular double-well potentials.  We  assume exact geometry and accurate numerical integration. The paper also assumes that for the given surface, Killing fields are known, while practical methods for the approximate Killing fields recovery were considered in \cite{OlshanskiiShakipov2026}. 

Finite element analysis of the surface Stokes problem is available for several fitted and unfitted methods, e.g.~\cite{OlshanskiiQuainiReuskenYushutin2018,fries2018higher,BonitoDemlowLicht2020,OlshanskiiReuskenZhiliakov2021,reusken2025analysis,hardering2025parametric,demlow2026taylor,kone2026divergence}. Structure-preserving methods for surface fluid--phase-field systems have also been developed~\cite{PalzhanovZhiliakovQuainiOlshanskii2021,BachiniKrauseVoigt2023,BachiniKrauseNitschkeVoigt2023}. These results do not by themselves give error estimates for the coupled  problem considered here. Related error analyses in Euclidean domains include a mixed finite element method for a Cahn--Hilliard--Darcy--Stokes system with a nonsteady flow equation~\cite{DiegelFengWise2015} and a SAV--MAC scheme for the Cahn--Hilliard--Stokes case~\cite{LiShen2020}. Our analysis addresses the surface system and, in particular, the nonlinear coupling and the Killing penalty.

The starting point is the continuous dependence estimate in energy norms proved in \cite{OlshanskiiShakipov2026}. Its discrete counterpart requires more than the usual $L^2$-based approximation theory: the nonlinear potential is controlled using mixed space--time norms involving $L^\infty(\G)$ and $W^{1,6}(\G)$. We therefore develop approximation estimates for TraceFEM in $L^p$- and $W^{1,p}$-based stabilized norms for $1\leq p\leq\infty$, and inverse estimates for $2\leq p\leq\infty$. The scaling of the normal derivative stabilization depends on $p$ and accounts for the codimension of the surface. An elliptic reconstruction then gives a discrete Sobolev embedding with a constant independent of the surface position in the mesh. These results extend the stabilized $L^2$ framework of~\cite{BouckNochettoShakipovYushutin2026} and are also useful for the analysis of other nonlinear surface problems.

Analysis of the paper requires a new discrete inf-sup condition, which applies to the velocity energy norm containing a Killing penalty with weight $h^{-2}$. We prove that the condition is equivalent to a more standard surface  inf-sup condition  for sufficiently small $h$. This transfers known stability results to the present method for the $P_{k+1}/P_k$ trace Taylor--Hood pairs, $k\geq1$, when $d=2$~\cite{OlshanskiiReuskenZhiliakov2021,JankuhnOlshanskiiReuskenZhiliakov2021}. For other pairs and for $d=3$, the analysis is conditional on the corresponding discrete inf-sup property.

Combining these tools with stabilized Ritz and Stokes projections yields a perturbation estimate for the difference between the Galerkin solution and the projected exact solution. Under the regularity assumptions stated in the error analysis and for balanced polynomial degrees, this difference converges faster in the energy norms than the total approximation error. We obtain optimal spatial orders for balanced polynomial degrees in the energy norms and in $L^\infty(0,T;L^2)$ for the phase field, $L^4(0,T;L^2)$ for the chemical potential, and $L^2(0,T;\bL^2)$ for the velocity. No additional duality argument for the nonlinear coupled problem is needed; duality is used in the projection estimates. For $d=2$, a discrete sharp Sobolev inequality also gives a pointwise phase-field error bound consisting of the pointwise Ritz-projection error and a higher-order term with a factor $|\log h|^{1/2}$. Numerical experiments examine the convergence rates for the fully discrete method with approximate geometry.

The remainder of the paper is organized as follows. Section~\ref{sec:surface-sch} recalls the model, its weak formulation, and the necessary analytical results. Section~\ref{sec:semidiscrete} presents the semidiscrete TraceFEM and its energy and conservation properties. Section~\ref{sec:tracefem} develops the approximation, projection, inverse, and inf-sup estimates used in the error analysis of Section~5. Section~6 reports numerical results. The appendix establishes an $L^p$ trace inequality and stability of the constant normal extension.

\section{The surface Stokes--Cahn--Hilliard problem}\label{sec:surface-sch}

Let $\G\subset\bbR^{d+1}$, $d=2,3$, be a closed, connected, sufficiently smooth hypersurface. 
Denote its unit normal by $\bn$ and the tangential projection by $\bP:=\bI-\bn\bn^T$. For a smooth scalar field $v$ and vector field $\bv$, the surface gradient and covariant derivative are
\[
    \nablaG v:=\bP\nabla\widetilde v,
    \qquad \nablaG\bv:=\bP(\nabla\widetilde{\bv})\bP,
    \qquad D_\G\bv:=\tfrac12\bigl(\nablaG\bv+(\nablaG\bv)^T\bigr),
\]
where the tildes denote any smooth extensions to a neighborhood of $\G$; the resulting surface operators do not depend on these extensions. We set $\divG\bv:=\tr(\nablaG\bv)$, apply surface divergence row-wise to tensors, and write $\DeltaG v:=\divG\nablaG v$ and $\Delta_S\bv:=\bP\divG D_\G\bv$~\cite{JankuhnOlshanskiiReusken2018}. Throughout, $(\cdot,\cdot)$ denotes the $L^2(\G)$ inner product, with Euclidean or Frobenius contraction for vector or tensor fields, and $\|\cdot\|:=\|\cdot\|_{L^2(\G)}$.

The inertia-free  surface Stokes--Cahn--Hilliard system reads: 
 Given the initial phase field $c_0$, the surface Stokes--Cahn--Hilliard problem is to find the tangential velocity $\bu$, pressure $p$, phase field $c$, and chemical potential $\mu$ satisfying
\begin{subequations}\label{eq:sch}
\begin{align}
    -2\eta\Delta_S\bu+\nablaG p&=-\sigma c\nablaG\mu,\label{eq:sch-mom}\\
    \divG\bu&=0,\label{eq:sch-div}\\
    \pa_t c+\bu\cdot\nablaG c-M\DeltaG\mu&=0,\label{eq:sch-phase}\\
    \mu&=\varepsilon^{-1}f_0'(c)-\varepsilon\DeltaG c,\label{eq:sch-chem}
\end{align}
\end{subequations}
with $c(0)=c_0$. The positive constants $\eta$, $M$, and $\varepsilon$ are the surface viscosity, mobility, and diffuse interface width, respectively; $\sigma>0$ is the line tension coefficient for $d=2$ (surface tension for $d=3$). The  double-well potential $f_0$ favors the two phases $c=\pm1$. We take $0<\varepsilon\leq1$ in the analysis.

The surface viscous operator has the finite-dimensional kernel consisting of infinitesimal rigid tangential motions of $\Gamma$, 
\[
E:=\{\bv \in H^1(\G)^{d+1} : \bv\cdot\bn=0,~~ D_\G \bv = 0 \text{ on } \G\},
\]
the space of Killing fields. The Killing fields are divergence-free. In \cite{OlshanskiiShakipov2026} we show that for sufficiently smooth $c$ and $\mu$ solving \eqref{eq:sch-chem},   the capillary force in \eqref{eq:sch-mom} is orthogonal to $E$,
	\begin{equation}\label{eq:orth}
	{(c\nablaG\mu,\bv_E)=0}\qquad  \text{for any }~ \bv_E\in E.
\end{equation}  
Hence the momentum equation is compatible  with this kernel and one may select a unique velocity by requiring $\bu\perp E$ in $\bL^2(\G)$.

For smooth solutions, testing the momentum, phase, and chemical potential equations with $\bu$, $\sigma\mu$, and $\sigma\pa_t c$, respectively, gives the energy balance: 
\begin{equation}\label{eq:en-law}
    \cE[c(t)]+\int_0^t\bigl[2\eta\|D_\G\bu(s)\|^2+\sigma M\|\nablaG\mu(s)\|^2\bigr]\,\mathrm{d}s=\cE[c_0],
    \qquad
    \cE[v]:=\sigma\left(\frac{\varepsilon}{2}\|\nablaG v\|^2+\varepsilon^{-1}(f_0(v),1)\right).
\end{equation}
Integrating \eqref{eq:sch-phase} over the closed surface also yields $(c(t),1)=(c_0,1)$. Thus the model dissipates interfacial free energy while conserving the total phase.

\smallskip
{We introduce the notation needed to formulate}~\eqref{eq:sch} weakly. Let $L^p(\G)$ and $W^{k,p}(\G)$ denote the usual surface Lebesgue and Sobolev spaces, and $H^k(\G):=W^{k,2}(\G)$. Boldface and double boldface indicate vector- and tensor-valued spaces. Subscripts $\#$ and $t$ indicate zero mean and tangential fields:
\begin{equation}\label{eq:subspaces}
    W^{k,p}_\#(\G):=\{v\in W^{k,p}(\G):(v,1)=0\},
    \qquad
    \bW^{k,p}_t(\G):=\{\bv\in\bW^{k,p}(\G):\bv\cdot\bn=0\}.
\end{equation}
We use the analogous notation $C^k_\#$ and $\bC^k_t$ for continuously differentiable functions, and $p'$ for the conjugate exponent. We shall use the Poincar\'e inequality
\begin{equation}\label{eq:poincare}
    \|v\|_{L^p(\G)}\leq C_P\|\nablaG v\|_{\bL^p(\G)},
    \qquad v\in W^{1,p}_\#(\G),\quad 1<p<\infty,
\end{equation}
with $C_P$ depending on $\G$ and $p$.

Consider the $L^2(\G)$-orthogonal decomposition  of $\bL^2_t(\G)$,
\begin{equation}\label{eq:decomp}
    \bL^2_t(\G)=E\oplus E^\perp,
    \qquad \cP_E:\bL^2_t(\G)\to E,
    \qquad \cP_E^\perp:=\operatorname{id}-\cP_E,
\end{equation}
where $\cP_E$ is the orthogonal projection. We recall the properties needed below from~\cite[Proposition~3.1]{OlshanskiiShakipov2026}.
\begin{proposition}[Killing fields]\label{prop:Killing}
    If $\G$ is of class $C^m$, $m\geq2$, then $E\subset\bW^{m-1,p}_t(\G)$ for every $1<p<\infty$. Moreover, the decomposition \eqref{eq:decomp} is stable in $\bH^{m-1}(\G)$:
    \begin{equation}\label{eq:stab-Killing-decomp}
        \|\cP_E\bv\|_{\bH^{m-1}(\G)}+\|\cP_E^\perp\bv\|_{\bH^{m-1}(\G)}\leq C\|\bv\|_{\bH^{m-1}(\G)},
        \qquad \bv\in\bH^{m-1}_t(\G),
    \end{equation}
 and $\cP_E$ extends to a bounded operator on $\bL^p_t(\G)$ for $1\leq p\leq\infty$.
\end{proposition}
Define
\begin{equation}\label{eq:bHk-tperp}
    \bW^{k,p}_{t,\perp}(\G):=\{\bv\in\bW^{k,p}_t(\G):\cP_E\bv=\bzero\},\qquad k\geq0.
\end{equation}
On $\bH^1_{t,\perp}(\G)$, the surface Korn inequality gives~\cite{JankuhnOlshanskiiReusken2018}
\begin{equation}\label{eq:korn}
    \|\bv\|_{\bH^1(\G)}\leq C_K\|D_\G\bv\|,
    \qquad \bv\in\bH^1_{t,\perp}(\G).
\end{equation}

Set $a(u,v):=(\nablaG u,\nablaG v)$ and $\ba(\bu,\bv):=(D_\G\bu,D_\G\bv)$. For a scalar function $w$, write $\mean{w}:=|\G|^{-1}(w,1)$ and $\zm{w}:=w-\mean{w}$. We define the \emph{centered} convection form 
\begin{equation}\label{eq:b-form}
    b(\bu,v,w):=(\bu\cdot\nablaG v,\zm{w}),
    \qquad (\bu,v,w)\in\bH^1(\G)\times H^1(\G)\times H^1(\G).
\end{equation}
It is continuous by $H^1(\G)\hookrightarrow L^4(\G)$ for $d\leq3$, and $b(\bu,v,1)=0$. If $\bu$ is tangential and divergence-free, it agrees with the usual convection form and is skew-symmetric in its last two arguments. 

Let $I:=(0,T)$ with $T>0$. Given $c_0\in H^1(\G)$, the weak problem is to find
\begin{equation}\label{eq:trialSpace}
\begin{alignedat}{2}
    &\bu\in L^2(I;\bH^1_{t,\perp}(\G)),\qquad &&c\in L^\infty(I;H^1(\G))\cap H^1(I;H^1(\G)'),\\
    &p\in L^2(I;L^2_\#(\G)),&&\mu\in L^2(I;H^1(\G)),
\end{alignedat}
\end{equation}
such that $c(0)=c_0$ and, for a.e. $t\in I$ and all $(\bw,q,v,g)\in\bH^1_{t,\perp}(\G)\times L^2_\#(\G)\times H^1(\G)\times H^1(\G)$,
\begin{subequations}\label{eq:weak}
\begin{align}
    2\eta \, \ba(\bu,\bw)-(p,\divG\bw)&=\sigma b(\bw,c,\mu),\label{eq:weak-mom}\\
    (q,\divG\bu)&=0,\label{eq:weak-div}\\
    \langle\pa_t c,v\rangle_{H^1(\G)'\times H^1(\G)}+b(\bu,c,v)+Ma(\mu,v)&=0,\label{eq:weak-phase}\\
    (\mu,g)&=\varepsilon^{-1}(f_0'(c),g)+\varepsilon a(c,g).\label{eq:weak-chem}
\end{align}
\end{subequations}
Note that using the centered trilinear form $\sigma b(\bw,c,\mu)$  for the capillary term  in the weak momentum equation amounts to the pressure shift $p\mapsto p+\sigma\zm{(c\zm{\mu})}$. We use $p$ for this shifted pressure below. This formulation preserves the cancellation of the coupling terms even when the discrete velocity is not pointwise divergence-free, and ensures unconditional phase conservation since $b(\bv, v, 1) = 0$ for all $(\bv, v) \in \bH^1_{t, \perp}(\G) \times H^1(\G)$.

The following results about the well-posedness of \eqref{eq:weak} and solution regularity were proved in \cite{OlshanskiiShakipov2026}:
\begin{proposition} [surface Stokes--Cahn--Hilliard] \label{prop:well-posedness}
Assume a  nonnegative double-well potential $f_0\in C^{2,1}_{loc}(\bbR)$ that satisfies
\begin{equation}\label{eq:double-well}
	|f_0'''(s)|\leq L(1+|s|)\quad\text{for a.e. }s\in\bbR 
\end{equation}
for some $L>0$.
 Assume  $\G$ is of class $C^2$, $d=2,3$. Then for every $c_0\in H^1(\G)$, problem \eqref{eq:weak} has a unique global weak solution. It satisfies \eqref{eq:en-law} and $(c(t),1)=(c_0,1)$ for $t\in[0,T]$. In addition,
    \[
        c\in L^2(I;W^{2,6}(\G))\cap L^4(I;L^\infty(\G)),\quad \text{and}\quad
        \qquad f_0'(c),\DeltaG c\in L^2(I;H^1(\G)).
    \]
  Two solutions with initial data $c_{0,1},c_{0,2}\in H^1(\G)$ satisfy the following continuous dependence bound
    \[
    \begin{split}
        &\|c_1-c_2\|_{L^\infty(I;H^1(\G))}+\|\bu_1-\bu_2\|_{L^2(I;\bH^1(\G))}
        +\|\mu_1-\mu_2\|_{L^2(I;H^1(\G))}\\
        &\qquad+\|p_1-p_2\|_{L^2(I;L^2(\G))}
        +\|\pa_t(c_1-c_2)\|_{L^2(I;H^1(\G)')}
        \leq C_T\|c_{0,1}-c_{0,2}\|_{H^1(\G)}.
    \end{split}
    \]
  Here $C_T$ depends on the surface, potential, material parameters, final time, and a priori bounds for the two solutions.
  Moreover,  if $d=2$, $m\geq0$ is an integer, $\G$ is of class $C^{m+4}$, $c_0\in H^{m+2}(\G)$, and $f_0\in C^{m+2,1}_{loc}(\bbR)$, then
    \[
    \begin{aligned}
        c&\in C([0,T];H^{m+2}(\G))\cap H^1(I;H^m(\G))\cap L^2(I;H^{m+4}(\G)),\\
        \mu&\in L^2(I;H^{m+2}(\G)),\qquad
        \bu\in L^2(I; \bH^{m+3}(\G)),\qquad p\in L^2(I; H^{m+2}(\G)).
    \end{aligned}
    \]
\end{proposition}
Note that the quartic potential $f_0(s)=\tfrac14(s^2-1)^2$ satisfies this assumption.

\section{Semidiscretization in space with TraceFEM}\label{sec:semidiscrete}

Let $\dst$ be the signed distance to $\G$, and choose $\delta_0>0$ such that the closest-point projection $\bp(\bx):=\bx-\dst(\bx)\bn(\bx)$ is well-defined in $U_{\delta_0}:=\{\bx:|\dst(\bx)|<\delta_0\}$. In this neighborhood $\bn=\nabla\dst$, and $v^e:=v\circ\bp$ denotes the constant normal extension, applied componentwise to vector and tensor fields. If $\G$ is of class $C^m$, $m\geq2$, then $\dst\in C^m(U_{\delta_0})$~\cite[Lemma~14.16]{GilbargTrudinger2001}. Let $\kappa_{\max}$ be the maximum absolute principal curvature, and fix
\begin{equation}\label{eq:assump-delta}
    0<\delta\leq\delta_0,\qquad 5\delta\kappa_{\max}<1.
\end{equation}
We use the normal and tangential projection extended to this tubular neighborhood.

Let $\{\widetilde\cT_h\}_h$ be a shape-regular, quasi-uniform family of simplicial meshes of a polyhedral domain $\Omega$ containing $\G$. Set $h:=\max_{T\in\widetilde\cT_h}\operatorname{diam}(T)$ and define
\begin{equation}\label{eq:cTh}
    \cT_h:=\{T\in\widetilde\cT_h:T\cap\G\neq\emptyset\},\qquad
    \Omega_h:=\bigcup_{T\in\cT_h}T.
\end{equation}
We take $h \leq \frac12 \min\{\delta,1\}$. In particular, $\Omega_h\subset U_\delta$, and the mesh resolves the curvature at the scale needed for the trace and extension estimates in the appendix. All surface integrals and geometric quantities in the analysis are exact. The notation $(\cdot,\cdot)_{\Omega_h}$ and $\|\cdot\|_{\Omega_h}$ refers to the bulk $L^2$ inner product and norm. We write $A\lesssim B$ or $A\simeq B$ for inequalities or equivalences with constants independent of $h$ and of the position of $\G$ in the mesh; other dependencies are stated when needed.

Let $\bV_h$, $Q_h$, and $V_h$ be the continuous bulk Lagrange spaces on $\cT_h$ of degrees $k_p+1$, $k_p$, and $k_c$, respectively, with $k_p,k_c\geq1$. Their traces are used to approximate velocity, pressure, and the two phase-field variables. We set $Q_{h,\#}:=Q_h\cap L^2_\#(\G)$ and $V_{h,\#}:=V_h\cap L^2_\#(\G)$.

For $l=0,1$, we define the following stabilization and penalty terms: 
\begin{equation} \label{eq:stab}
	\begin{split}
		& s_l^{\mathrm{ext}}(u , v ) := h^{1-2l} (\bn \cdot \nabla u , \bn \cdot \nabla v )_{\Omega_h}, \\
		& \bs_l^{\mathrm{ext}}(\bu , \bv ) := h^{1-2l} ((\nabla \bu ) \bn, (\nabla \bv ) \bn)_{\Omega_h}, \\
		& \bs^{\mathrm{tan}}(\bu , \bv ) := h^{-2} (\bn \cdot \bu , \bn \cdot \bv )_{\G},
	\end{split}
\end{equation}
for all $u,v \in H^1(\Omega_h)$, and all $\bu , \bv  \in \bH^1(\Omega_h)$.
To control the Killing modes of the discrete velocity, we {include} the additional Killing penalty term:
\begin{equation}
	\bs^{\mathrm{Kil}}(\bu , \bv ) := h^{-2} (\cP_E  \bu , \cP_E  \bv )_{\G}.
\end{equation}
The projection $\cP_E$ here is understood as an extension to $\bL^2(\G)$ by $\cP_E=(\cP_E)|_{\bL^2_t(\G)}\bP$. {The same extension convention is used for} $\cP_E^\perp$. With $a$ and $\ba$ from Section~\ref{sec:surface-sch} we define
\begin{equation}\label{eq:bilinear-forms}
\begin{aligned}
    (u,v)_*&:=(u,v)+s_0^{ext}(u,v),\\
    a_*(u,v)&:=a(u,v)+s_1^{ext}(u,v),\\
    \ba_*(\bu,\bv)&:=\ba(\bP\bu,\bP\bv)+\bs^{\mathrm{tan}}(\bu,\bv)
        +\bs_1^{ext}(\bu,\bv)+\bs^{\mathrm{Kil}}(\bu,\bv).
\end{aligned}
\end{equation}
The associated scalar mass and energy norms and the velocity energy norm are
\[
    \|v \|_{L^2_*}^2:=(v ,v )_*,\qquad
    |v |_{H^1_*}^2:=a_*(v ,v ),\qquad
    \|v \|_{H^1_*}^2:=\|v \|_{L^2_*}^2+|v |_{H^1_*}^2,
    \qquad \|\bv \|_{\bH^1_*}^2:=\ba_*(\bv ,\bv ).
\]
Section~\ref{sec:tracefem} recalls the bulk norm control and extends these definitions to the $L^p$ setting as needed for the analysis. The stabilizations vanish on normal extensions of exact scalar fields and of tangential velocity fields orthogonal to $E$, so the forms are consistent with \eqref{eq:weak}.

Given $c_0\in H^1(\G)$, the semidiscrete FE method is to find
\[
    (\bu_h,p_h,c_h,\mu_h)\in L^2(I;\bV_h)\times L^2(I;Q_{h,\#})\times H^1(I;V_h)\times L^2(I;V_h)
\]
such that $c_h(0)=R_hc_0$, with elliptic projection $R_h$ defined in \eqref{eq:h1-proj}, and
\begin{equation}\label{eq:semidiscrete}
\begin{split}
    2\eta \, \ba_*(\bu_h,\bw_h)+(\nablaG p_h,\bw_h)
        &=\sigma b(\cP_E^\perp\bw_h,c_h,\mu_h),\\
    (\nablaG q_h,\bu_h)-(2\eta)^{-1}s_0^{ext}(p_h,q_h)&=0,\\
    (\pa_t c_h,v_h)_*+b(\cP_E^\perp\bu_h,c_h,v_h)+Ma_*(\mu_h,v_h)&=0,\\
    (\mu_h,g_h)_*&=\varepsilon^{-1}(f_0'(c_h),g_h)+\varepsilon a_*(c_h,g_h),
\end{split}
\end{equation}
for a.e. $t\in I$ and all $(\bw_h,q_h,v_h,g_h)\in\bV_h\times Q_{h,\#}\times V_h\times V_h$\rev{.}
Note that since $b$ uses a tangential gradient in its second argument, its value is unchanged if the transport velocity $\cP_E^\perp\bw_h$ is replaced by $\bw_h-\cP_E\bw_h$. The factor $(2\eta)^{-1}$ in the pressure stabilization matches the viscosity scaling of the Stokes block.

For inf-sup stable Stokes trace finite elements, cf. Assumption~\ref{assump:disc-inf-sup}, and sufficiently small $h$, the FE method has a unique global solution~\cite{OlshanskiiShakipov2026}. 
Furthermore, the finite element method is structure preserving. Indeed, testing \eqref{eq:semidiscrete} with $\bw_h=\bu_h$, $q_h=p_h$, $v_h=\sigma\mu_h$, and $g_h=\sigma\pa_t c_h$ yields the discrete analogue of the energy balance \eqref{eq:en-law},
\begin{equation}\label{eq:discr-en-law}
    \ddt\cE_*[c_h(t)]+2\eta\|\bu_h\|_{\bH^1_*}^2
    +(2\eta)^{-1}s_0^{ext}(p_h,p_h)+\sigma Ma_*(\mu_h,\mu_h)=0,
\end{equation}
where
\[
    \cE_*[v_h]:=\sigma\left(\frac{\varepsilon}{2}a_*(v_h,v_h)+\varepsilon^{-1}(f_0(v_h),1)\right);
\]
{the centered convection form also gives conservation of the total phase} without any compatibility condition between $\bV_h$ and $V_h$
\begin{equation}\label{eq:discr-mass-conserv}
    \ddt(c_h,1)=-b(\cP_E^\perp\bu_h,c_h,1)=0.
\end{equation}

\section{Properties of stabilized TraceFEM}\label{sec:tracefem}

We collect tools needed for the error analysis. 

\subsection{Stabilized norms}

We first recall the bulk norm control provided by normal derivative stabilization. For either scalar space $W_h\in\{V_h,Q_h\}$,
\begin{equation}\label{eq:bulk-norm-control}
    \|v_h\|_{L^2_*}^2\simeq h^{-1}\|v_h\|_{\Omega_h}^2,\qquad v_h\in W_h.
\end{equation}
This follows from~\cite[Lemmas~7.8, 7.9]{GrandeLehrenfeldReusken2018} and the bulk inverse and trace inequalities. In particular, $(\cdot,\cdot)_*$ is an inner product on $W_h$. The Poincar\'e inequality \eqref{eq:poincare} carries over to the stabilized norms on $V_h$: for $v_h$ with vanishing mean over $\G$, it bounds the surface part of $\|v_h\|_{L^2_*}$ by $\|\nablaG v_h\|$, while the identity $s_0^{ext}=h\,s_1^{ext}$ and $h\leq1$ bound the extension part by $s_1^{ext}(v_h)$. Consequently, we obtain
\begin{equation}\label{eq:discr-poincare}
    \|v_h\|_{L^2_*}\lesssim|v_h|_{H^1_*}\quad\forall\,v_h\in V_{h,\#},
    \qquad
    \|v_h\|_{H^1_*}\lesssim|v_h|_{H^1_*}+|\mean{v_h}|\quad\forall\,v_h\in V_h,
\end{equation}
with constants depending only on $C_P$ and $|\G|$.

For the nonlinear analysis we extend this framework to $L^p$-based quantities. With $1\leq p\leq\infty$ 
define
\begin{equation}\label{eq:stab-seminorm}
\begin{aligned}
    s_{l,p}^{ext}(v)&:=h^{1-l-1/p}\|\bn\cdot\nabla v\|_{L^p(\Omega_h)}, & l &=-1,0,1,\\
    \bs_{l,p}^{ext}(\bv)&:=h^{1-l-1/p}\|(\nabla\bv)\bn\|_{\bL^p(\Omega_h)}, & l&=0,1,\\
    \bs_{l,p}^{tan}(\bv)&:=h^{-l}\|\bn\cdot\bv\|_{L^p(\G)}, \quad
    \bs_{l,p}^{Kil}(\bv):=h^{-l}\|\cP_E\bv\|_{\bL^p(\G)},\quad & l & =0,1.
\end{aligned}
\end{equation}
To be consistent with the earlier notation of stabilization, we will suppress  $p$ when $p=2$; for the bilinear stabilizations defined in \eqref{eq:stab}, the corresponding seminorm is their square root on the diagonal. 
The factor $h^{-1/p}$ accounts for the change from bulk to surface measure. In particular, H\"older's inequality gives
\[
    |s_0^{ext}(u,v)|\leq s_{l,p}^{ext}(u)s_{-l,p'}^{ext}(v),\qquad l=0,1.
\]
For finite $p$, set
\begin{equation}\label{eq:*-norms}
\begin{aligned}
    \|v_h\|_{L^p_*}&:=\bigl(\|v_h\|_{L^p(\G)}^p+s_{0,p}^{ext}(v_h)^p\bigr)^{1/p},\\
    |v_h|_{W^{1,p}_*}&:=\bigl(\|\nablaG v_h\|_{\bL^p(\G)}^p+s_{1,p}^{ext}(v_h)^p\bigr)^{1/p},\\
    \|v_h\|_{W^{1,p}_*}&:=\bigl(\|v_h\|_{L^p_*}^p+|v_h|_{W^{1,p}_*}^p\bigr)^{1/p},\\
    \|v_h\|_{W^{-1,p}_*}&:=\bigl(\|v_h|_\G\|_{W^{1,p'}(\G)'}^p+s_{-1,p}^{ext}(v_h)^p\bigr)^{1/p}.
\end{aligned}
\end{equation}
For $p=\infty$, use the sum of the two component quantities in each definition. We write $H^j_*=W^{j,2}_*$, $j=-1,0,1$, with $H^0_*=L^2_*$. In particular, the negative surface norm is the dual of the full $H^1(\G)$ space, including constants. These definitions imply
\[
    |(u_h,v_h)_*|\leq\|u_h\|_{W^{l,p}_*}\|v_h\|_{W^{-l,p'}_*},\qquad l=0,1.
\]
The vector quantities needed below are
\begin{equation}\label{eq:*-norms-l2}
\begin{aligned}
    \|\bv_h\|_{\bL^p_*}&:=\bigl(\|\bv_h\|_{\bL^p(\G)}^p+\bs_{0,p}^{ext}(\bv_h)^p\bigr)^{1/p},\\
    \|\bv_h\|_{\bW^{1,p}_*}&:=\bigl(\|D_\G(\bP\bv_h)\|_{\dL^p(\G)}^p
        +\bs_{1,p}^{tan}(\bv_h)^p+\bs_{1,p}^{ext}(\bv_h)^p+\bs_{1,p}^{Kil}(\bv_h)^p\bigr)^{1/p},
\end{aligned}
\end{equation}
again using sums when $p=\infty$. The case $p=2$ agrees with the velocity energy norm in Section~\ref{sec:semidiscrete}. Orthogonality gives the useful identity
\begin{equation}\label{eq:L2*-norm-expanded}
    \|\bv_h\|_{\bL^2_*}^2
    =\|\cP_E^\perp\bv_h\|^2+\bs_0^{Kil}(\bv_h)^2
      +\bs_0^{tan}(\bv_h)^2+\bs_0^{ext}(\bv_h)^2.
\end{equation}
\begin{lemma}
	\label{lem:norms-bVh}
    The quantities $\|\cdot\|_{\bL^2_*}$ and $\|\cdot\|_{\bH^1_*}$ are norms on $\bV_h$, and
    \[
        \|\bv_h\|_{\bL^2_*}+\|\cP_E^\perp\bv_h\|_{\bH^1(\G)}
        \lesssim\|\bv_h\|_{\bH^1_*}.
    \]
    The scalar and vector quantities in \eqref{eq:*-norms} and \eqref{eq:*-norms-l2} are also norms on their bulk finite element spaces for $1\leq p\leq\infty$. The scalar seminorm $|\cdot|_{W^{1,p}_*}$ is a norm on the zero-mean spaces for $1<p<\infty$.
\end{lemma}
\begin{proof}
    That $\|\cdot\|_{\bL^2_*}$ defines a norm follows componentwise from \eqref{eq:bulk-norm-control}. Korn's inequality applied to $\cP_E^\perp\bv_h$, together with \eqref{eq:L2*-norm-expanded} and $\bs_0^X=h\bs_1^X$ for $X\in\{ext,tan,Kil\}$, gives the asserted estimate and hence $\|\cdot\|_{\bH^1_*}$ is a norm. For general $p$, vanishing of a scalar norm means that both the surface trace  and the bulk normal derivative vanish. Equation \eqref{eq:bulk-norm-control} then gives $v_h=0$. In the vector energy case, vanishing strain and Killing penalty imply $\bP\bv_h=0$; the tangential penalty gives $\bn\cdot\bv_h=0$, and the bulk mass control completes the argument. Poincar\'e's inequality handles the scalar zero-mean subspaces.
\end{proof}

{In stabilized forms and error norms, surface functions are understood as their constant normal extensions. Since their normal derivatives vanish, the scalar stabilized norms coincide with the corresponding surface norms. For $\bv\in\bW^{1,p}_{t,\perp}(\G)$, the tangential and Killing penalties also vanish, so $\|\bv\|_{\bW^{1,p}_*}=\|D_\G\bv\|_{\dL^p(\G)}$.}

Before introducing projections, define the discrete dual norm
\begin{equation}\label{eq:discr-dual-norm}
    \|v_h\|_{V_h'}:=\sup_{0\neq\phi_h\in V_h}
        \frac{|(v_h,\phi_h)_*|}{\|\phi_h\|_{H^1_*}}.
\end{equation}
This is the norm of the functional induced by the stabilized mass pairing.

\subsection{Interpolation, projection and inverse estimates}

Let $I_h:W^{1,1}(\Omega_h)\to W_h$ be the Scott--Zhang \cite{ScottZhang1990} quasi-interpolation operator into a scalar Lagrange space of degree $\ell\geq1$, and let $\bI_h$ act componentwise. We use its standard local stability and approximation properties. The mesh restriction in Section~\ref{sec:semidiscrete} ensures that Lemmas~\ref{lem:local-trace-ineq} and \ref{lem:const-normal-ext} apply to the active elements and the extension neighborhoods used below.

\begin{lemma}[interpolation] \label{lem:approx}
            If $\G$ is of class $C^{m+1}$ for an integer $1 \leq m \leq \ell + 1$, and $u \in W^{m, p}(\G), p \in [1,\infty]$, then the following estimate is valid:
            \begin{equation} \label{eq:approx-scal}
                \|u - I_h u^e\|_{L^p_*} + h \|u - I_h u^e\|_{W^{1,p}_*} \leq C h^m \|u\|_{W^{m,p}(\G)}.
            \end{equation}
            Under the same assumptions on $m$ and $\G$, if $\bu \in \bW^{m, p}_t(\G)$, then the following estimates are valid:
            \begin{equation} \label{eq:approx-vec}
                \begin{split}
                    \|\bu - \bI_h \bu^e\|_{\bL^p_*} + h\, \big(\|\nabla_M (\bu - \bP \bI_h \bu^e)\|_{\dL^p(\G)} + \bs_{1,p}^{tan}(\bI_h \bu^e) + \bs_{1,p}^{ext}(\bI_h \bu^e)\big) & \leq C h^m \|\bu\|_{\bW^{m, p}(\G)},
                \end{split}
            \end{equation}
            where $\nabla_M \bu$ is componentwise $\nablaG$. If $\bu \in \bW^{m, p}_{t, \perp}(\G)$, then
            \begin{equation} \label{eq:approx-Kil}
                \|\bu - \bI_h \bu^e\|_{\bW^{1,p}_*} \leq C h^{m-1} \|\bu\|_{\bW^{m, p}(\G)}.
            \end{equation}
        \end{lemma}
        \begin{proof}
    Let $\theta=u^e-I_hu^e$. For $m\geq2$, the local trace inequality from Lemma~\ref{lem:local-trace-ineq}, applied to $\theta$ and its elementwise gradient, and the identity $s_{0,p}^{ext}=h s_{1,p}^{ext}$ give
    \[
    \begin{split}
        \|\theta\|_{L^p_*}+h\|\theta\|_{W^{1,p}_*}
        &\lesssim h^{-1/p}\|\theta\|_{L^p(\Omega_h)}
        +h^{1-1/p}\|\nabla\theta\|_{\bL^p(\Omega_h)}
        +h^{2-1/p}\|\nabla_h\nabla\theta\|_{\dL^p(\Omega_h)}\\
        &\lesssim h^{m-1/p}\|u^e\|_{W^{m,p}(U_{2h})}
        \lesssim h^m\|u\|_{W^{m,p}(\G)}.
    \end{split}
    \]
    Here $\nabla_h$ is the elementwise gradient, and the last two steps use Scott--Zhang approximation and Lemma~\ref{lem:const-normal-ext}. When $m=1$, the scalar trace error is estimated in the same way using $W^{1,p}$ stability. For the surface gradient, one instead uses
    \[
        \|\nablaG\theta\|_{\bL^p(\G)}
        \leq\|\nablaG u\|_{\bL^p(\G)}
        +C h^{-1/p}\|\nabla I_hu^e\|_{\bL^p(\Omega_h)}
        \lesssim\|u\|_{W^{1,p}(\G)},
    \]
    by Corollary~\ref{cor:local-trace-discrete} and extension stability. The stabilization seminorm is bounded directly by $h^{1-1/p}\|\nabla\theta\|_{\bL^p(\Omega_h)}$. These arguments also hold at $p=\infty$.

    Componentwise application gives the vector mass and stabilization estimates. Tangentiality gives
    \[
        \bs_{1,p}^{tan}(\bI_h\bu^e)
        \leq h^{-1}\|\bu-\bI_h\bu^e\|_{\bL^p(\G)}.
    \]
    Since $\bu-\bP\bI_h\bu^e=\bP(\bu-\bI_h\bu^e)$, the product rule and boundedness of $\nabla_M\bP$ give the componentwise surface gradient estimate in \eqref{eq:approx-vec}. Finally, if $\cP_E\bu=0$, Proposition~\ref{prop:Killing} implies
    \[
        \bs_{1,p}^{Kil}(\bI_h\bu^e)
        =h^{-1}\|\cP_E(\bI_h\bu^e-\bu)\|_{\bL^p(\G)}
        \lesssim h^{-1}\|\bI_h\bu^e-\bu\|_{\bL^p(\G)}.
    \]
    This proves \eqref{eq:approx-Kil}.
\end{proof}

        Let $P_h : L^2(\G) \to V_h$ be the stabilized $L^2$-projection, namely, given $u \in L^2(\G)$, find $P_h u \in V_h$ such that 
        \begin{equation} \label{eq:L2-proj}
            (P_h u, v_h)_* = (u, v_h)\quad \forall\,v_h \in V_h.
        \end{equation}
        The following result was proven in \cite[Lemma~2.6 and Eq.~(5.15)]{BouckNochettoShakipovYushutin2026}.
        \begin{proposition} \label{prop:Ph}
            $P_h$ defined in \eqref{eq:L2-proj} satisfies the following properties:
            \begin{enumerate}
                \item $L^2_*$-stability: $\|P_h\|_{\cL(L^2_*)} := \|P_h\|_{\cL(L^2(\G), (V_h, L^2_*))} \leq 1$.
                \item $H^1_*$-stability: $\|P_h\|_{\cL(H^1_*)} := \|P_h\|_{\cL(H^1(\G), (V_h, H^1_*))} \leq C$, where $C > 0$ depends only on $\G$ and shape-regularity, but not on $h$ or how $\G$ cuts the bulk mesh.
            \end{enumerate}
        \end{proposition}
        
    \medskip
        For all $u, v \in H^1(\G)$ and all $u_h, v_h \in V_h$, we define
        \begin{equation} \label{eq:h1-prod}
            (u, v)_{H^1(\G)} := (u, v) + a(u, v), \qquad (u_h, v_h)_{H^1_*} := (u_h, v_h)_* + a_*(u_h, v_h).
        \end{equation}
        Let $R_h : H^1(\G) \to V_h$ be the stabilized $H^1$-projection, namely, given $u \in H^1(\G)$, find $R_h u \in V_h$ such that for all $v_h \in V_h$,
        \begin{equation} \label{eq:h1-proj}
            (R_h u, v_h)_{H^1_*} = (u, v_h)_{H^1(\G)}.
        \end{equation}
        The Ritz projection has the following standard properties.
        \begin{lemma} \label{lem:prop-Rh}
            $R_h$ satisfies the following properties:
            \begin{enumerate}
                \item Mean preservation: $(R_h u, 1) = (u, 1)$ for all $u \in H^1(\G)$.
                \item Projection identity: $(R_h u - u, v_h)_* = a_*(u - R_h u, v_h)$ for all $v_h \in V_h$.
                \item $H^1_*$-stability: $\|R_h\|_{\cL(H^1(\G), (V_h, H^1_*))} \leq 1$.
                \item Error estimate: for $\G$ of class $C^{\max\{m+1,3\}}, u \in H^m(\G)$, and integer $1 \leq m \leq k_c + 1$, the following error estimate is valid:
                \begin{equation} \label{eq:h1-proj-err}
                    \|u - R_h u\|_{L^2_*} + h \|u - R_h u\|_{H^1_*} \leq C h^m \|u\|_{H^m(\G)}.
                \end{equation}
            \end{enumerate}
        \end{lemma}
        \begin{proof}
            Mean preservation follows by testing with $1$, and the next two properties follow directly from \eqref{eq:h1-proj}. Galerkin orthogonality and Lemma~\ref{lem:approx} give the $H^1_*$ error bound. The surface $L^2$ bound follows by the usual duality argument for $I-\DeltaG$ and $H^2$ elliptic regularity. Finally, $s_0^{ext}=h s_1^{ext}$ bounds the remaining stabilization term.
        \end{proof}

        The following result is a generalization of \cite[Lemma 5.8]{BouckNochettoShakipovYushutin2026} to the $L^p$-setting for $p \in [2, \infty]$ that will be needed in the error analysis.
        \begin{lemma}[inverse estimates] \label{lem:stab-inverse}
            There exists $C > 0$ depending only on $\G, p, k_c$ and shape regularity such that for $k = 0,1$, $p \in [2, \infty]$, and all $v_h \in V_h$,
            \begin{equation} 
                \begin{split}
                    \|v_h\|_{W^{k,p}_*} & \leq C h^{- d/2 - (k - d/p)} \|v_h\|_{L^2_*}.
                \end{split}
            \end{equation}
        \end{lemma}
\begin{proof}
    We set $1/\infty:=0$ and use $h_T\simeq h$ and $h\leq1$. Corollary~\ref{cor:local-trace-discrete}, applied on each $T\in\cT_h$ to $v_h$ and to $\nabla v_h$, together with $|\nablaG v_h|\leq|\nabla v_h|$, bounds the surface terms in $\|v_h\|_{W^{k,p}_*}$ by $h^{-1/p}\|v_h\|_{W^{k,p}(\Omega_h)}$. The stabilization terms satisfy $s_{l,p}^{ext}(v_h)\leq h^{1-l-1/p}\|\nabla v_h\|_{\bL^p(\Omega_h)}$ with $l\leq k$. Bulk inverse estimates in dimension $d+1$ and \eqref{eq:bulk-norm-control} then give
    \[
        \|v_h\|_{W^{k,p}_*}
        \lesssim h^{-k-1/p}\|v_h\|_{L^p(\Omega_h)}
        \lesssim h^{-k-1/p-(d+1)(1/2-1/p)}\|v_h\|_{\Omega_h}
        \lesssim h^{-k-d(1/2-1/p)}\|v_h\|_{L^2_*}.
    \]
    This is the assertion, since $-k-d(1/2-1/p)=-d/2-(k-d/p)$.
\end{proof}

    \subsection{Discrete Sobolev inequality}
        For $v_h \in V_h$ we let $\cA_h v_h \in V_h$ be the discrete $\cA := I - \DeltaG$ operator, given by
        \begin{equation} \label{eq:cAh}
            (\cA_h v_h, w_h)_* = (v_h, w_h)_{H^1_*}, \qquad \forall w_h \in V_h.
        \end{equation}
        Using $\cA_h$, we can define the discrete $H^{2}(\G)$-norm:
        \begin{equation} \label{eq:H2*}
            \|v_h\|_{H^2_*} := \|\cA_h v_h\|_{L^2_*}, \qquad \forall v_h \in V_h.
        \end{equation}
        We need a FE analogue of the embedding $H^2(\G)\hookrightarrow W^{1,6}(\G)$. For this,  we use a continuous elliptic reconstruction. For $v_h \in V_h$, the \textit{elliptic reconstruction} $v = \cR v_h \in H^1(\G)$ is given by
        \begin{equation} \label{eq:cR}
            (v, w)_{H^1(\G)} = (\cA_h v_h, w), \quad w \in H^1(\G).
        \end{equation}
        In other words, $\cR v_h = \cA^{-1} \cA_h v_h$. In flat domains, the analog of \eqref{eq:cR} was used in \cite{HeywoodRannacher1982} to prove a discrete Agmon inequality, while \cite{MakridakisNochetto2003} used it for a posteriori error estimates for parabolic problems. Below we derive its properties in the TraceFEM setting, noting that $R_h \cR v_h \neq v_h$ because of the stabilization. 
        \begin{lemma} \label{lem:elliptic-reconstruction}
            Let $\G$ be of class $C^3$, and $v_h \in V_h$ be given. The elliptic reconstruction $v = \cR v_h \in H^1(\G)$ from \eqref{eq:cR} satisfies $v \in H^2(\G)$ and 
            \begin{equation} \label{eq:ell-reconstr}
                \|v - v_h\|_{L^2_*} \leq C h^2 \|v_h\|_{H^2_*}, \qquad \|v\|_{H^2(\G)} \leq C \|v_h\|_{H^2_*},
            \end{equation}
            where $C > 0$ depends only on $\G$, shape-regularity and $\deg V_h$.
        \end{lemma}
        \begin{proof}
            By elliptic regularity we have $\|v\|_{H^2(\G)} \leq C \|\cA_h v_h\| \leq C \|v_h\|_{H^2_*}$. Note that the stability constant is independent of $h$ or how $\G$ cuts the bulk mesh. For the $L^2_*$-error estimate, we add and subtract $R_h v$ and use the $L^2_*$-error estimate for it followed by $\|v\|_{H^2(\G)} \leq C\|v_h\|_{H^2_*}$:
            \begin{equation} \label{eq:ell-rec-l2-est-step}
                \|v^e - v_h\|_{L^2_*} \leq C h^2 \|v\|_{H^2(\G)} + \|R_h v - v_h\|_{L^2_*} \leq C h^2 \|v_h\|_{H^2_*} + \|R_h v - v_h\|_{L^2_*}.
            \end{equation}
            Note that the last quantity is not zero. Testing \eqref{eq:cR} with $w = w_h \in V_h$ and using the definition \eqref{eq:h1-proj} of $R_h$, we have
            \begin{equation} \label{eq:cAhvh-vs-Rhv}
                (\cA_h v_h, w_h) = (v, w_h)_{H^1(\G)} = (R_h v, w_h)_{H^1_*}, \quad \forall w_h \in V_h.
            \end{equation}
            Thus, by \eqref{eq:cAh} and \eqref{eq:cAhvh-vs-Rhv} the defect $e_h = v_h - R_h v \in V_h$ satisfies
            \begin{equation} \label{eq:agmon-eh}
                (e_h, w_h)_{H^1_*} = (\cA_h v_h, w_h)_* - (\cA_h v_h, w_h) = s_0^{ext}(\cA_h v_h, w_h), \quad \forall w_h \in V_h.
            \end{equation}
            We claim that $\|e_h\|_{L^2_*} \leq C h^2 \|v_h\|_{H^2_*}$. Testing \eqref{eq:agmon-eh} with $w_h = e_h$ and using $s_0^{ext}(e_h) = h \, s_1^{ext}(e_h) \leq h \|e_h\|_{H^1_*}$ gives $\|e_h\|_{H^1_*} \leq h \|\cA_h v_h\|_{L^2_*} = h \|v_h\|_{H^2_*}$, and therefore $s_0^{ext}(e_h) \leq C h^2 \|v_h\|_{H^2_*}$. For the part on $\G$ we argue by duality: let $\phi \in H^1(\G)$ solve $(\phi, w)_{H^1(\G)} = (e_h, w)$ for all $w \in H^1(\G)$, so that $\|\phi\|_{H^2(\G)} \leq C \|e_h\|$. Then, by \eqref{eq:h1-proj}, \eqref{eq:agmon-eh}, the fact that $s_0^{ext}(\cA_h v_h, \phi^e) = 0$ because $\bn \cdot \nabla \phi^e = 0$, and Lemma \ref{lem:prop-Rh},
            \[
                \begin{split}
                    \|e_h\|^2 & = (\phi, e_h)_{H^1(\G)} = (R_h \phi, e_h)_{H^1_*} = s_0^{ext}(\cA_h v_h, R_h \phi) = s_0^{ext}(\cA_h v_h, R_h \phi - \phi^e) \\
                    & \leq \|\cA_h v_h\|_{L^2_*} \, h \, \|\phi - R_h \phi\|_{H^1_*} \leq C h^2 \|v_h\|_{H^2_*} \|e_h\|.
                \end{split}
            \]
            Therefore, $\|e_h\|_{L^2_*} \leq C h^2 \|v_h\|_{H^2_*}$, which combined with \eqref{eq:ell-rec-l2-est-step} finishes the proof.
        \end{proof}
        We next use the elliptic reconstruction to derive a discrete Sobolev inequality.
        \begin{lemma}[discrete Sobolev embedding] \label{lem:disc-sob}
            Let $\G$ be of class $C^3$. Then there exists a constant $C > 0$ depending only on $\G$, shape-regularity and $\deg V_h$ such that for all $v_h \in V_h$, $\|v_h\|_{W^{1,6}_*} \leq C \|v_h\|_{H^2_*}$.
        \end{lemma}
        \begin{proof}
            Let $v_h \in V_h$ and $v \in H^2(\G)$ be its elliptic reconstruction from \eqref{eq:cR}. We add and subtract $I_h v^e$, use its stability in $W^{1,6}_*$ and stabilized inverse estimates (Lemma \ref{lem:stab-inverse}) as follows:
            \[
                \begin{split}
                    \|v_h\|_{W^{1,6}_*} & \leq \|I_h v^e\|_{W^{1,6}_*} + \|v_h - I_h v^e\|_{W^{1,6}_*} \lesssim \|v\|_{W^{1,6}(\G)} + h^{-1-d/3} (\|v_h - v\|_{L^2_*} + \|I_h v^e - v\|_{L^2_*}) \lesssim \|v_h\|_{H^2_*},
                \end{split}
            \]
            where we used error estimates for $I_h$, properties of elliptic reconstruction \eqref{eq:ell-reconstr}, $d \leq 3$, and Sobolev embedding. This finishes the proof.
        \end{proof}
        The next important ingredient is the following pointwise bound on the stabilized $H^1$ projection. 
        \begin{lemma}[$L^\infty$-control of $R_h$] \label{lem:Linfty-Rh}
            Let $\G$ be of class $C^3$. Then there exists $C > 0$, depending only on $\G$, $k_c$ and shape-regularity, such that for all $v \in H^2(\G)$,
            \[
                \|R_h v\|_{L^\infty(\Omega_h)} \leq C \|v\|_{H^1(\G)}^{1/2} \|v\|_{H^2(\G)}^{1/2}.
            \]
        \end{lemma}
        \begin{proof}
            Let $w_h := R_h v - I_h v^e \in V_h$. Lemma \ref{lem:prop-Rh}  and Lemma \ref{lem:approx}  give 
            \[
                \|w_h\|_{L^2_*} \leq C \min\left\{h \|v\|_{H^1(\G)}, h^2 \|v\|_{H^2(\G)}\right\} \leq C h^{3/2} \|v\|_{H^1(\G)}^{1/2} \|v\|_{H^2(\G)}^{1/2}.
            \]
            The bulk inverse estimate and $\|\cdot\|_{L^2(\Omega_h)} \leq C h^{1/2} \|\cdot\|_{L^2_*}$ \cite[Lemmas 7.8, 7.9]{GrandeLehrenfeldReusken2018} yield $\|w_h\|_{L^\infty(\Omega_h)} \leq C h^{-d/2} \|w_h\|_{L^2_*} \leq C h^{(3-d)/2} \|v\|_{H^1(\G)}^{1/2} \|v\|_{H^2(\G)}^{1/2}$. Here, Agmon's inequality on $\G$, $\|v\|_{L^\infty(\G)} \leq C \|v\|_{H^1(\G)}^{1/2} \|v\|_{H^2(\G)}^{1/2}$ for $d \leq 3$, follows from the classical one on $\bbR^d$ by localization with a finite $C^2$ atlas and a subordinate partition of unity. Finally, the $L^\infty(\Omega_h)$-stability of the Scott--Zhang operator, Lemma \ref{lem:const-normal-ext} ($W^{k,p}$-stability of the constant normal extension), and Agmon's inequality on $\G$ give $\|I_h v^e\|_{L^\infty(\Omega_h)} \leq C \|v\|_{L^\infty(\G)} \leq C \|v\|_{H^1(\G)}^{1/2} \|v\|_{H^2(\G)}^{1/2}$. The claim follows from $d \leq 3$ and $h \leq 1$.
        \end{proof}
        The following result is from \cite[Lemma~4.1]{PalzhanovZhiliakovQuainiOlshanskii2021}. It was stated for smooth surfaces in $\bbR^3$, though only $C^2$ regularity was actually used. The proof in fact controls $\|v_h\|_{L^\infty(\Omega_h)}$ rather than $\|v_h\|_{L^\infty(\G)}$; see Eq. (4.7)-(4.8) therein. 
        \begin{lemma} \label{lem:discr-sharp-Sobolev}
            Let $\G$ be of class $C^2$ and $d=2$. Then there exists $C > 0$ depending only on $\G$, $\deg V_h$ and shape-regularity such that for all $v_h \in V_h$,
            \begin{equation}
                \|v_h\|_{L^\infty(\Omega_h)} \leq C |\log h|^{1/2} \|v_h\|_{H^1_*}.
            \end{equation}
        \end{lemma}

    \subsection{Inf-sup condition and Stokes projection} \label{subsec:disc-inf-sup}
        We are interested in a suitable inf-sup stability property of the velocity--pressure pair $(\bV_h, Q_{h, \#})$. We introduce the velocity norm (note that the Killing stabilization $\bs_1^{Kil}(\cdot,\cdot)$ is not included)
        \begin{equation} \label{eq:tilde-bH1-norm}
            \|\bw_h\|_{\wt \bH^1_*} := \left( \|\bP \bw_h\|_{\bH^1(\G)}^2 + \bs_1^{tan}(\bw_h)^2 + \bs_1^{ext}(\bw_h)^2 \right)^{1/2}, \qquad \bw_h \in \bV_h,
        \end{equation}
       and assume that the velocity and pressure spaces satisfy the following condition:
        \begin{assump}\label{assump:disc-inf-sup}
            There exist $\alpha_0, \beta_0 > 0$, independent of $h > 0$ and of how the hypersurface cuts the ambient mesh, such that for all $q_h \in Q_{h, \#}$,
            \begin{equation} \label{eq:discr-inf-sup-ORZ}
                \alpha_0 \|q_h\|_{L^2_*} \leq \sup_{\bzero \neq \bw_h \in \bV_h} \frac{(\nablaG q_h, \bw_h)}{\|\bw_h\|_{\wt \bH^1_*}} + \beta_0 \, s_0^{ext}(q_h).
            \end{equation}
        \end{assump}
        
        For {$d=2$ and sufficiently small $h$}, the condition \eqref{eq:discr-inf-sup-ORZ} is known to hold for the $P_{k+1}/P_k$, $k\ge1$ trace elements~\cite{OlshanskiiReuskenZhiliakov2021,JankuhnOlshanskiiReuskenZhiliakov2021}, i.e. generalized trace Taylor--Hood elements.
        
        We also note that the stabilization term on the right-hand side of \eqref{eq:discr-inf-sup-ORZ} is 
        necessary for \eqref{eq:discr-inf-sup-ORZ} to hold as seen from the following example: Consider the unit sphere and a pressure space of degree at least two, with $q_h(\bx) := |\bx|^2 - 1$ which vanishes identically on $\G$ but not in $\Omega_h$. Thus, in this example $\nablaG q_h = \bzero$ but $\|q_h\|_{L^2_*}^2 = s_0^{ext}(q_h)^2 > 0$.

        For the stability of our finite element method we need, however, a stronger inf-sup condition: 
        \begin{equation} \label{eq:discr-inf-sup}
            \alpha \|q_h\|_{L^2_*} \leq \sup_{\bzero \neq \bw_h \in \bV_h} \frac{(\nablaG q_h, \bw_h)}{\|\bw_h\|_{\bH^1_*}} + \beta \, s_0^{ext}(q_h), \qquad \forall q_h \in Q_{h, \#}.
        \end{equation}
        with the $\|\cdot\|_{\bH^1_*}$ norm in the denominator. 
        
        Korn's inequality and finite dimensionality of $E$ imply
        \[
            \|\bP\bw_h\|_{\bH^1(\G)}^2
            \simeq\|D_\G(\bP\bw_h)\|^2+\|\cP_E\bw_h\|^2.
        \]
        Thus the two velocity norms differ only in the weight of the Killing component: $1$ in $\|\cdot\|_{\wt \bH^1_*}$ and $h^{-2}$ in $\|\cdot\|_{\bH^1_*}$. The next lemma shows that \eqref{eq:discr-inf-sup-ORZ} yields \eqref{eq:discr-inf-sup}.

        \begin{lemma} \label{lem:discr-inf-sup}
            Let $\G$ be of class $C^3$. Let $(\bV_h, Q_{h, \#})$ be any finite element pair with $\ell := \deg \bV_h \geq 1$ satisfying Assumption \ref{assump:disc-inf-sup}. Then there exist $h_0 > 0$ depending only on $\alpha_0, \G, d, \ell$ and shape-regularity, and $\alpha, \beta > 0$ depending only on $\alpha_0, \beta_0, \G, d, \ell$ and shape-regularity, such that \eqref{eq:discr-inf-sup} holds for all $h \leq h_0$. 
        \end{lemma}
\begin{proof}
    Let $S(q_h)$ denote the supremum on the right-hand side of \eqref{eq:discr-inf-sup}. For $\bv_h\in\bV_h$ with $\|\bv_h\|_{\wt\bH^1_*}=1$, set $\bk:=\cP_E\bv_h$ and $\bw_h:=\bv_h-\bI_h\bk^e$. Proposition~\ref{prop:Killing}, finite dimensionality of $E$, and $\G\in C^3$ imply $\|\bk\|_{\bH^2(\G)}\lesssim1$. Since $\divG\bk=0$, integration by parts and \eqref{eq:approx-vec} with $m=2$ give
    \begin{equation}\label{eq:nablaG-qh-bIh-bk}
        |(\nablaG q_h,\bI_h\bk^e)|
        =|(q_h,\divG(\bP\bI_h\bk^e-\bk))|
        \lesssim h\|q_h\|_{L^2_*}.
    \end{equation}
    The first three terms in $\|\bw_h\|_{\bH^1_*}$ are bounded by the triangle inequality and interpolation stability. For the Killing term, $\cP_E\bv_h=\bk$ and \eqref{eq:approx-vec} yield
    \[
        \bs_1^{Kil}(\bw_h)
        =h^{-1}\|\cP_E(\bk-\bI_h\bk^e)\|
        \lesssim h^{-1}\|\bk-\bI_h\bk^e\|
        \lesssim h\|\bk\|_{\bH^2(\G)}.
    \]
    Thus $\|\bw_h\|_{\bH^1_*}\lesssim1$, and, whether or not $\bw_h=0$,
    \[
        |(\nablaG q_h,\bv_h)|
        \leq |(\nablaG q_h,\bw_h)|+|(\nablaG q_h,\bI_h\bk^e)|
        \lesssim S(q_h)+h\|q_h\|_{L^2_*}.
    \]
    Taking the supremum over $\bv_h$ and using \eqref{eq:discr-inf-sup-ORZ} proves \eqref{eq:discr-inf-sup} after absorbing the last term for sufficiently small $h$. Conversely, Korn's inequality gives $\|\bv_h\|_{\wt\bH^1_*}\lesssim\|\bv_h\|_{\bH^1_*}$ for $h\leq1$, so \eqref{eq:discr-inf-sup} implies \eqref{eq:discr-inf-sup-ORZ}.
\end{proof}

        Let $(\bS_h, S_h) : \bH^1(\G) \times L^2(\G) \to \bV_h \times Q_{h, \#}$ be the stabilized Stokes projection, namely, given $(\bu, p) \in \bH^1(\G) \times L^2(\G)$, find $(\bS_h \bu, S_h p) \in \bV_h \times Q_{h, \#}$ such that for all $(\bv_h, q_h) \in \bV_h \times Q_{h, \#}$,
        \begin{equation} \label{eq:Stokes-proj}
            \begin{split}
                2 \eta \, \ba_*(\bS_h \bu, \bv_h) + (\nablaG S_h p, \bv_h) & = 2 \eta \, \ba(\bu, \bP \bv_h) - (p, \divG (\bP \bv_h)), \\
                (\nablaG q_h, \bS_h \bu) - (2 \eta)^{-1} s_0^{ext}(S_h p, q_h) & = (\nablaG q_h, \bu).
            \end{split}
        \end{equation}
        Here, $\bS_h \bu$ and $S_h p$ both depend on $(\bu, p)$; this abuse of notation will be used for brevity.
        \begin{lemma} \label{lem:prop-Stokes-proj}
            Under Assumption~\ref{assump:disc-inf-sup} and for $h\leq h_0$ from Lemma~\ref{lem:discr-inf-sup}, the stabilized Stokes projection defined in \eqref{eq:Stokes-proj} satisfies the following properties:
            \begin{enumerate}
                \item \begingroup
Stability: there is $C>0$, independent of $h$ and of how $\G$ intersects the bulk mesh, such that
                \begin{equation}\label{eq:Stokes-proj-stability}
                \|\bS_h\bu\|_{\bH^1_*}+\|S_hp\|_{L^2_*}
                \leq C\bigl(\|\bu\|_{\bH^1(\G)}+\|p\|_{L^2(\G)}\bigr).
                \end{equation}
                The constant may depend on $\eta$ and the discrete inf-sup constants.
\endgroup
                \item Error estimates: for $\G$ of class $C^{m+2}$, $\bu\in \bH^1_{t, \perp}(\G)\cap \bH^{m+1}(\G)$, $ p \in L^2_\#\cap H^m(\G)$, and integer $1 \leq m \leq k_p+1$, the following approximation estimate is valid:
                \begin{equation} \label{eq:Stokes-err}
                    \|\bu - \bS_h \bu\|_{\bL^2_*} + h(\|\bu - \bS_h \bu\|_{\bH^1_*} + \|p - S_h p\|_{L^2_*}) \leq C h^{m+1} (\|\bu\|_{\bH^{m+1}(\G)} + \|p\|_{H^m(\G)}).
                \end{equation}
            \end{enumerate}
        \end{lemma}
        \begin{proof}
            The stability estimate \eqref{eq:Stokes-proj-stability} follows from coercivity of $\ba_*$, the discrete inf-sup condition \eqref{eq:discr-inf-sup}, and continuity of the right-hand side in \eqref{eq:Stokes-proj}. The defining equations give Galerkin orthogonality:
            \begin{equation} \label{eq:stokes-proj-galerkin-orth}
                \begin{split}
                    2 \eta \, \ba_*(\bS_h \bu - \bu, \bv_h) - (S_h p - p, \divG (\bP \bv_h)) & = 0, \\
                    - (q_h, \divG (\bP \bS_h \bu - \bu)) - (2 \eta)^{-1} s_0^{ext}(S_h p, q_h) & = 0.
                \end{split}
            \end{equation}
            Applying uniform mixed stability to the difference between the Stokes projection and arbitrary finite element approximants, and using these identities, yields
            \[
            \|\bu-\bS_h\bu\|_{\bH^1_*}+\|p-S_hp\|_{L^2_*}
            \lesssim \inf_{\bv_h\in\bV_h}\|\bu-\bv_h\|_{\bH^1_*}
            +\inf_{q_h\in Q_{h,\#}}\|p-q_h\|_{L^2_*}.
            \]
            Choosing Scott--Zhang approximants, with the pressure mean subtracted, and applying Lemma~\ref{lem:approx} gives the asserted energy estimate. Next we apply the duality argument: consider $(\bphi, \psi) \in \bH^1_{t, \perp}(\G) \times L^2_\#(\G)$ solving $2 \eta \, \ba(\bphi, \bv) - (\psi, \divG \bv) + (\divG \bphi, q) = (\bu - \cP_E^\perp \bS_h \bu, \bv)$ for all $(\bv, q) \in \bH^1_{t, \perp}(\G) \times L^2_\#(\G)$. For $\G$ of class $C^3$, higher Stokes regularity \cite[Theorem 6.1]{BenavidesNochettoShakipov2025-vector} ensures that $\|\bphi\|_{\bH^2(\G)} + \|\psi\|_{H^1(\G)} \leq C(\G,\eta^{-1}) \|\bu - \cP_E^\perp \bS_h \bu\|$. Then using $D_\G \bv_E = \bzero$ and $\divG \bv_E = 0$ for all $\bv_E \in E$, we infer
            \[
                \begin{split}
                    \|\bu - \cP_E^\perp \bS_h \bu\|^2 & = |2 \eta \, \ba(\bphi, \bu - \cP_E^\perp \bS_h \bu) - (\psi, \divG (\bu - \cP_E^\perp \bS_h \bu))| \\
                    & = |2 \eta \, \ba_*(\bphi, \bu - \bS_h \bu) - (\psi, \divG(\bu - \bP \bS_h \bu))|.
                \end{split}
            \]
            We next use Galerkin orthogonality \eqref{eq:stokes-proj-galerkin-orth} together with $\divG \bphi = 0$ to write 
            \[
                \begin{split}
                    |2 \eta \, & \ba_*(\bphi, \bu - \bS_h \bu) - (\psi, \divG(\bu - \bP \bS_h \bu))| \\
                    & = |2 \eta \, \ba_*(\bphi - \bS_h \bphi, \bu - \bS_h \bu) - (S_h p - p, \divG (\bP \bS_h \bphi - \bphi)) \\
                    & \quad - (\psi - S_h \psi, \divG(\bu - \bP \bS_h \bu)) - (2 \eta)^{-1} s_0^{ext}(S_h p, S_h \psi)| \\
                    & \leq 2 \eta \|\bphi - \bS_h \bphi\|_{\bH^1_*} \|\bu - \bS_h \bu\|_{\bH^1_*} + \sqrt{d} \|S_h p - p\|  \|\bphi - \bS_h \bphi\|_{\bH^1_*} \\
                    & \quad + \sqrt{d} \|\psi - S_h \psi\| \|\bu - \bS_h \bu\|_{\bH^1_*} + (2 \eta)^{-1} s_0^{ext}(S_h p) s_0^{ext}(S_h \psi) \\
                    & \leq \sqrt2 (\sqrt{d} + 2 \eta + (2 \eta)^{-1}) (\|\bphi - \bS_h \bphi\|_{\bH^1_*} + \|\psi - S_h \psi\|_{L^2_*}) (\|\bu - \bS_h \bu\|_{\bH^1_*} + \|p - S_h p\|_{L^2_*}).
                \end{split}
            \]
            Using $\bH^1_* \times L^2_*$-error estimate from \eqref{eq:Stokes-err} together with $\|\bphi\|_{\bH^2(\G)} + \|\psi\|_{H^1(\G)} \leq C(\G,\eta^{-1}) \|\bu - \cP_E^\perp \bS_h \bu\|$ yields $\|\bu - \cP_E^\perp \bS_h \bu\| \leq C h^{m+1} (\|\bu\|_{\bH^{m+1}(\G)} + \|p\|_{H^m(\G)})$. We next note that $\bs_0^{ext}(\bS_h \bu) + \bs_0^{tan} (\bS_h \bu) + \bs_0^{Kil}(\bS_h \bu)$ is $O(h^{m+1})$ due to the $\bH^1_*$ error estimate in \eqref{eq:Stokes-err} and the scaling relation $\bs_0^{X}(\cdot) = h \bs_1^X(\cdot), X \in \{ext, tan, Kil\}$. The $\bL^2_*$-error estimate then follows from identity \eqref{eq:L2*-norm-expanded}. This finishes the proof.
        \end{proof}

\section{Error analysis}
    In this section, we carry out a numerical analysis of the semidiscrete method \eqref{eq:semidiscrete}. First, we define a suitable interpolant of the exact solution using Stokes and $H^1$-projections introduced earlier. This interpolant satisfies the Galerkin scheme up to residual terms (defects). We then adapt the proof of the Lipschitz stability bound from \cite[Thm.~4.4]{OlshanskiiShakipov2026} to the discrete setting. This amounts to deriving discrete regularity estimates for the Galerkin solution, which are then used to prove a discrete perturbation estimate. Estimating the defects then lets us derive error estimates in the energy norm. Using interpolation inequalities and Sobolev embeddings, we also derive estimates in weaker norms.

    All assumptions from the previous section remain in force.  Also throughout this section, assumption \eqref{eq:double-well} on the double-well potential plays an important role. 
We will need  the following consequences of \eqref{eq:double-well}:
\begin{equation}\label{eq:f0-consequences}
	\begin{aligned}
		|f_0''(s)|&\leq L'(1+s^2),&
		|f_0''(a)-f_0''(b)|&\leq L'(1+|a|+|b|)|a-b|,\\
		|f_0'(s)|&\leq L'(1+|s|^3),&
		|f_0'(a)-f_0'(b)|&\leq L'(1+a^2+b^2)|a-b|,
	\end{aligned}
\end{equation}
Here $L'$ depends only on $L$, $|f_0'(0)|$, and $|f_0''(0)|$.

    \subsection{Solution interpolants and some bounds}
        We define $(\wbu_h, \wp_h) := (\bS_h \bu, S_h p) \in \bV_h \times Q_{h, \#}$ using the Stokes projection from \eqref{eq:Stokes-proj}, and $(\wc_h, \wmu_h) := (R_h c, R_h \mu) \in V_h \times V_h$, using the $H^1$-projection from \eqref{eq:h1-proj}. Testing the momentum equation in \eqref{eq:weak} with $\bw = \cP_E^\perp \bw_h \in \bH^1_{t, \perp}(\G)$, and noting that $b(\bw_h - \cP_E \bw_h, v, g) = b(\cP_E^\perp \bw_h, v, g)$ for all $v, g \in H^1(\G)$ because the second argument enters through $\nablaG$, we find that $(\wbu_h, \wp_h, \wc_h, \wmu_h)$ satisfy for a.e. $t \in I$ and all $\bw_h \in \bV_h$,
        \begin{equation} \label{eq:sD-bu}
            \begin{split}
                2 \eta \, \ba_*(\wbu_h, \bw_h) + (\nablaG \wp_h, \bw_h) & - \sigma b(\cP_E^\perp \bw_h, \wc_h, \wmu_h) \\
                & = \sigma b(\cP_E^\perp \bw_h, c, \mu) - \sigma b(\cP_E^\perp \bw_h, \wc_h, \wmu_h) =: \sD_\bu(\bw_h).
            \end{split}
        \end{equation}
        Using property (2) of $R_h$ from Lemma \ref{lem:prop-Rh}, we find that $(\wc_h, \wmu_h)$ give rise to the following defects:
        \begin{equation} \label{eq:sD-c-mu}
            \begin{split}
                (\pa_t \wc_h, v_h)_* & + b(\cP_E^\perp \wbu_h, \wc_h, v_h) + M a_*(\wmu_h, v_h) \\
                & = (\pa_t \wc_h, v_h)_* - \langle \pa_t c, v_h \rangle + b(\cP_E^\perp \wbu_h, \wc_h, v_h) - b(\bu, c, v_h) + M (\mu - \wmu_h, v_h)_* =: \sD_c(v_h), \\
                (\wmu_h, g_h)_* & - \varepsilon^{-1} (f_0'(\wc_h), g_h) - \varepsilon a_*(\wc_h, g_h) \\
                & = (\wmu_h - \mu, g_h)_* + \varepsilon^{-1} (f_0'(c) - f_0'(\wc_h), g_h) - \varepsilon (c - \wc_h, g_h)_* =: \sD_\mu(g_h).
            \end{split}
        \end{equation}
         We note that the divergence equation carries no defect since $(\nablaG q_h, \bu) = - (q_h, \divG \bu)$ for $\bu \in \bH^1_t(\G)$. Finally, for a fixed $g_h \in V_h$, we compute
        \begin{equation} \label{eq:pat-sD-mu}
            \begin{split}
                \pa_t \sD_\mu (g_h) & = (\pa_t [\wmu_h - \mu], g_h)_* + \varepsilon^{-1} (f''_0(c) \pa_t c - f''_0(\wc_h) \pa_t \wc_h, g_h) - \varepsilon (\pa_t [c - \wc_h], g_h)_* \\
                & = (\pa_t [\wmu_h - \mu], g_h)_* + \varepsilon^{-1} (f''_0(c) \pa_t (c - \wc_h) + [f''_0(c) - f''_0(\wc_h)] \pa_t \wc_h, g_h) - \varepsilon (\pa_t [c - \wc_h], g_h)_*.
            \end{split}
        \end{equation}
        
        \begin{remark}
            Using the stabilized $H^1$-projection allows converting $a_*(\wc_h - c, \cdot)$ into $(c - \wc_h, \cdot)_*$ and $a_*(\wmu_h - \mu, \cdot)$ into $(\mu - \wmu_h, \cdot)_*$. Hence, no $\|\nablaG (c - \wc_h)\|$ appears in $\sD_\mu$ and no $\|\nablaG (\mu - \wmu_h)\|$ appears in $\sD_c$.
        \end{remark}
        \begin{remark}
            The defects are estimated in Lemma \ref{lem:defects} and require Assumption \ref{assump:disc-inf-sup} ($(\bV_h, Q_{h, \#})$ is discretely inf-sup stable) with $h \leq h_0$ as well as additional regularity of the exact solution.
        \end{remark}
    
        For every $\bv_h \in \bV_h$, we define $\bv_{h,\perp} = \cP_E^\perp \bv_h \in \bH^1_{t, \perp}(\G)$.
        \begin{remark}[$\bH^1(\G)$-control of $\bv_{h,\perp}$] \label{rmk:bH1-control}
            Let $\bv \in \bH^1_{t, \perp}(\G)$. Then, Korn's inequality \eqref{eq:korn} yields
            \[
            \begin{split}
                \|\bv_{h, \perp} - \bv\|_{\bH^1(\G)} & \leq C_K \|D_\G (\bv_{h, \perp} - \bv)\| = C_K \|D_\G (\cP_E^\perp [\bv_h - \bv])\| = C_K \|D_\G (\bP [\bv_h - \bv])\| \leq C_K \|\bv_h - \bv\|_{\bH^1_*}.
            \end{split}
            \]
        \end{remark}
        \begin{lemma}[continuity of $b$] \label{lem:b-disc}
            For all $\bz_h \in \bV_h$ and all $v, w \in H^1(\G)$,
            \begin{equation} \label{eq:b-disc}
                |b(\bz_{h, \perp}, v, w)| \leq C_b \|\bz_h\|_{\bH^1_*} \|\nablaG v\| \|\nablaG w\|.
            \end{equation}
        \end{lemma}
        \begin{proof}
            By H\"older's inequality we get
            $
                |b(\bz_{h,\perp}, v, w)| \leq \|\bz_{h,\perp}\|_{\bL^4(\G)} \|\nablaG v\| \|\zm{w}\|_{L^4(\G)}.
            $
            Remark \ref{rmk:bH1-control} with $\bv = \bzero$ and the Sobolev embedding give $\|\bz_{h, \perp}\|_{\bL^4(\G)} \leq \Csob C_K \|\bz_h\|_{\bH^1_*}$, while Poincar\'e \eqref{eq:poincare} gives $\|\zm{w}\|_{L^4(\G)} \leq \Csob (1 + C_P) \|\nablaG w\|$.
        \end{proof}
        The following result is the foundation for the perturbation error estimate.
        \begin{proposition}[discrete a priori bounds] \label{prop:discrete-regularity}
            Let $\G$ be of class $C^3$. The finite element solution $\{\bu_h, c_h, \mu_h\}$ and the interpolants $\{\wt \bu_h,\wt c_h,\wt \mu_h\}$ satisfy the bounds
            \begin{equation} \label{eq:discrete-regularity}
            	\begin{split}
\|\bu_h\|_{L^2(I; \bH^1_*)} + \|c_h\|_{L^\infty(I; H^1_*)} + \|\mu_h\|_{L^2(I; H^1_*)} +~ &\|c_h\|_{L^2(I; W^{1,6}_*)} \leq C,\\
\|\wt \bu_h\|_{L^2(I; \bH^1_*)} + \|\wt c_h\|_{L^\infty(I; H^1_*)} + \|\wt \mu_h\|_{L^2(I; H^1_*)} +~ & \|\wt c_h\|_{L^4(I; L^\infty(\Omega_h))} \leq C,
            	\end{split}
             \end{equation}
            where $C$ is independent of $h$ and depends only on $\G, T, \sigma, M, \varepsilon^{-1}, \eta^{-1}, L'$, $(c_0, 1)$, norms of the exact solution from Proposition \ref{prop:well-posedness} (surface Stokes--Cahn--Hilliard), and discrete inf-sup constants from \eqref{eq:discr-inf-sup}.
        \end{proposition}
        Before proceeding with the proof, we note {the different norms} for $c_h$ and $\wt c_h = R_h c$ in \eqref{eq:discrete-regularity}. These norms are sufficient to control the nonlinearity of $f_0$ (cf. Lemma \ref{lem:ophih} and the subsequent discussion).
        \begin{proof}
            Recall from Proposition \ref{prop:well-posedness} that $\bu \in L^2(I; \bH^1(\G)), c \in L^2(I; H^2(\G)) \cap L^4(I; L^\infty(\G)) \cap L^\infty(I; H^1(\G))$ and $\mu \in L^2(I; H^1(\G))$.
            
            The bounds on $\wbu_h, \wc_h$ and $\wmu_h$ in $H^1_*$ follow from $H^1$-stability of the respective projections. The bounds on $\wc_h$ in $L^4(I; L^\infty(\Omega_h))$ follow from Lemma \ref{lem:Linfty-Rh} ($L^\infty$-control of $R_h$). We now address the Galerkin solution. The discrete energy law \eqref{eq:discr-en-law} yields control, uniform in $h$, of 
            \[
                \sup_{t \in I} |c_h(t)|_{H^1_*}^2 + \int_I (\|\bu_h\|_{\bH^1_*}^2 + |\mu_h|_{H^1_*}^2 ) \ds.
            \]
            The control of $\|\bu_h\|_{L^2(I; \bH^1_*)}$ is immediate. The mean of $c_h$ is controlled by $\mean{c_0}$ due to the conservation of total phase, which, by \eqref{eq:discr-poincare}, gives control of $\|c_h\|_{L^\infty(I; H^1_*)}$. For $\mu_h$, we test the chemical potential equation in \eqref{eq:semidiscrete} with $g_h = 1$ and use \eqref{eq:f0-consequences} to deduce
            \[
                |(\mu_h, 1)| \leq \varepsilon^{-1} |(f'_0(c_h), 1)| \leq \varepsilon^{-1} |\G|^{1/2} \|f_0'(c_h)\| \leq C(1 + \|c_h\|_{L^6(\G)}^3) \leq C (1 + \|c_h\|_{L^\infty(I; H^1(\G))}^3).
            \]
            This and \eqref{eq:discr-poincare} yield uniform control of $\|\mu_h\|_{L^2(I; H^1_*)}$. 
            
            {It remains to bound} $\|c_h\|_{L^2(I; W^{1,6}_*)}$. To that end, we need to establish a uniform bound on $\|c_h\|_{L^2(I; H^2_*)} = \|\cA_h c_h\|_{L^2(I; L^2_*)}$. By the last equation of \eqref{eq:semidiscrete} and the definition \eqref{eq:L2-proj} of $P_h$, $\varepsilon a_*(c_h, g_h) = (\mu_h - \varepsilon^{-1} P_h f_0'(c_h), g_h)_*$ for all $g_h \in V_h$. Therefore,
            \[
                \cA_h c_h = c_h + \varepsilon^{-1} \mu_h - \varepsilon^{-2} P_h f_0'(c_h), \qquad \|\cA_h c_h\|_{L^2_*} \leq \|c_h\|_{L^2_*} + \varepsilon^{-1} \|\mu_h\|_{L^2_*} + \varepsilon^{-2} \|f_0'(c_h)\|,
            \]
            where we used $\|P_h\|_{\cL(L^2(\G), (V_h, L^2_*))} \leq 1$. Note that all the terms on the right-hand side have already been estimated, which yields the asserted uniform control of $\|c_h\|_{L^2(I; H^2_*)}$. The bound on $\|c_h\|_{L^2(I; W^{1,6}_*)}$ then follows from Lemma \ref{lem:disc-sob}. This finishes the proof.
        \end{proof}
    
    \subsection{A discrete error estimate}
        Let $(\wbu_h, \wp_h, \wc_h, \wmu_h)$ be given as in the previous subsection and $(\obu_h, \op_h, \oc_h, \omu_h) := (\bu_h - \wbu_h, p_h - \wp_h, c_h - \wc_h, \mu_h - \wmu_h)$ be the discrete errors. Then these errors satisfy
        \begin{equation} \label{eq:err-eq}
            \begin{split}
                2 \eta \, \ba_*(\obu_h, \bw_h) + (\nablaG \op_h, \bw_h) & = \sigma b(\bw_{h, \perp}, c_h, \mu_h) - \sigma b(\bw_{h, \perp}, \wc_h, \wmu_h) - \sD_\bu(\bw_h), \\
                (\nablaG q_h, \obu_h) - (2 \eta)^{-1} s_0^{ext}(\op_h, q_h) & = 0, \\
                (\pa_t \oc_h, v_h)_* + M a_*(\omu_h, v_h) & = - b(\bu_{h, \perp}, c_h, v_h) + b(\wbu_{h, \perp}, \wc_h, v_h) - \sD_c(v_h), \\
                (\omu_h, g_h)_* - \varepsilon a_*(\oc_h, g_h) & = \varepsilon^{-1} (\ophi_h, g_h)_* - \sD_\mu(g_h),
            \end{split}
        \end{equation}
        where we used the notation $\bv_{h,\perp} := \cP_E^\perp \bv_h$ for $\bv_h \in \bV_h$, and 
        \begin{equation} \label{eq:ophi}
            \ophi := f_0'(c_h) - f_0'(\wc_h) \notin V_h, \qquad \ophi_h := P_h \ophi \in V_h.
        \end{equation}
        We note that conservation of mass together with $c_h(0) = R_h c_0$ implies $\oc_h \in V_{h,\#} := V_h \cap L^2_\#(\G)$ and $\oc_h(0) = 0$. 

        The key to the continuous dependence result at the PDE level \cite[Thm.~4.4]{OlshanskiiShakipov2026} was testing the phase-field error equation with $f_0'(c_1) - f_0'(c_2)$. {Such a test function is not admissible} at the discrete level, since $\ophi \notin V_h$. Instead we test with $\ophi_h = P_h \ophi$. The following result is the discrete counterpart of \cite[Lem.~4.3]{OlshanskiiShakipov2026}.
        \begin{lemma}[control of the nonlinearity] \label{lem:ophih}
            Let $\G$ be of class $C^3$. Then for a.e. $t \in I$,
            \[
                \|\ophi_h\|_{L^2_*}^2 \leq C_{\phi,0} |\oc_h|_{H^1(\G)}^2, \qquad \int_0^t \|\ophi_h\|_{H^1_*}^2 \ds \leq \int_0^t C_{\phi, 1}(s) |\oc_h|_{H^1(\G)}^2 \ds, 
            \]
            where $C_{\phi, 0} > 0$ depends on $\G, f_0$, and the $L^\infty(I; H^1(\G))$-norms of $c_h, \wc_h$, while 
            \begin{equation} \label{eq:Cphi1}
                C_{\phi, 1}(t) = C' \left[1 + \|\wc_h(t)\|_{L^\infty(\G)}^4 + (1 + \|c_h\|_{L^\infty(I; H^1(\G))}^2 + \|\wc_h\|_{L^\infty(I; H^1(\G))}^2) \|c_h(t)\|_{W^{1,6}(\G)}^2 \right]
            \end{equation}
            belongs to $L^1(I)$ with some $C' > 0$ depending only on $\G, f_0$, $\deg V_h$ and shape-regularity.
        \end{lemma}
        The asymmetry in the norms for $c_h$ and $\wc_h$ in \eqref{eq:Cphi1} comes directly from the PDE level. At the discrete level, this asymmetry lets us avoid a discrete Agmon inequality to control  $c_h$ in $L^4(I; L^\infty(\Omega_h))$.
        \begin{proof}
            The result for $\ophi_h$ replaced by $\ophi$ follows directly from \cite[Lemma~4.3]{OlshanskiiShakipov2026} and Proposition \ref{prop:discrete-regularity}. The result for $\ophi_h$ then follows from the $L^2_*$ and $H^1_*$-stability of $P_h$; cf. Proposition \ref{prop:Ph}. 
        \end{proof}

        For generic functionals $F, \bF$ on $V_h, \bV_h$, respectively, we define
        \begin{equation} \label{eq:discr-dual-norms}
            \|F\|_{V_h'} := \sup_{0 \neq \phi_h \in V_h} \frac{F(\phi_h)}{\|\phi_h\|_{H^1_*}}, \qquad \|\bF\|_{\bV_h'} := \sup_{\bzero \neq \boldsymbol{\phi}_h \in \bV_h} \frac{\bF( \boldsymbol{\phi}_h)}{\|\boldsymbol{\phi}_h\|_{\bH^1_*}}.
        \end{equation}
        This will be used for defects below, which brings us to the first main result of this work.
        
        \begin{theorem}[discrete error estimate] \label{thm:disc-energy}
            Let $\G$ be of class $C^3$ and $c, \mu \in H^1(I; H^1(\G))$. Then
            \begin{equation} \label{eq:disc-energy-est}
                |\oc_h(t)|_{H^1_*}^2 + \int_0^t \|\obu_h\|_{\bH^1_*}^2 + s_0^{ext}(\op_h)^2 + \|\omu_h\|_{H^1_*}^2 \ds \leq C_T^2 \cE_\sD^2,
            \end{equation}
            holds for every $t \in [0,T]$, with
            \begin{equation}
                 \cE_\sD := \left(\|\sD_\mu\|_{H^1(I; V_h')}^2 + \|\sD_c\|_{L^2(I; V_h')}^2 + \|\sD_\bu\|_{L^2(I; \bV_h')}^2 \right)^{1/2},
            \end{equation}
            and $C_T > 0$ is bounded uniformly in $h$ and depends only on $T, \G$, $\eta, \eta^{-1}, \sigma, M, M^{-1}, \varepsilon^{-1}$, $L'$, mesh shape regularity, norms on $(\bu, p, c, \mu)$ from Proposition \ref{prop:well-posedness} (surface Stokes--Cahn--Hilliard) and the constant from Proposition \ref{prop:discrete-regularity}.
        \end{theorem}
        \begin{proof}
            We test the error equation \eqref{eq:err-eq} with $\bw_h = \obu_h, q_h = \op_h, v_h = \sigma \omu_h, g_h = \sigma \pa_t \oc_h$ to obtain 
            \begin{equation*}
                \begin{split}
                    2 \eta \|\obu_h\|_{\bH^1_*}^2 + (2 \eta)^{-1} s_0^{ext}(\op_h)^2 & = \sigma b(\obu_{h, \perp}, c_h, \omu_h) + \sigma b(\obu_{h, \perp}, \oc_h, \wmu_h) - \sD_\bu(\obu_h), \\
                    \sigma (\pa_t \oc_h, \omu_h)_* + \sigma M |\omu_h|_{H^1_*}^2 & = - \sigma b(\obu_{h, \perp}, c_h, \omu_h) - \sigma b(\wbu_{h, \perp}, \oc_h, \omu_h) - \sigma \sD_c(\omu_h), \\
                    \sigma (\omu_h, \pa_t \oc_h)_* - \sigma \varepsilon a_*(\oc_h, \pa_t \oc_h) & = \sigma \varepsilon^{-1} (\ophi_h, \pa_t \oc_h)_* - \sigma \sD_\mu(\pa_t \oc_h),
                \end{split}
            \end{equation*}
            Combining the three equations and  the cancellation of $\sigma b(\obu_{h, \perp}, c_h, \omu_h)$ and $a_*(\oc_h, \pa_t \oc_h) = \frac12 \ddt |\oc_h|_{H^1_*}^2$ yields
            \begin{equation*}
                \begin{split}
                    \frac{\sigma \varepsilon}{2} & \ddt |\oc_h|_{H^1_*}^2 + 2 \eta \|\obu_h\|_{\bH^1_*}^2 + (2 \eta)^{-1} s_0^{ext}(\op_h)^2 + \sigma M |\omu_h|_{H^1_*}^2 \\
                    & = \sigma b(\obu_{h, \perp}, \oc_h, \wmu_h) - \sigma b(\wbu_{h, \perp}, \oc_h, \omu_h) - \sigma \varepsilon^{-1} (\ophi_h, \pa_t \oc_h)_* - \sD_\bu(\obu_h) - \sigma \sD_c(\omu_h) + \sigma \sD_\mu(\pa_t \oc_h).
                \end{split}
            \end{equation*}
            Integrating from $0$ to $t$ and using $\oc_h(0) = 0$ results in
            \begin{equation} \label{eq:err-eqn-strong}
                \begin{split}
                    \frac{\sigma \varepsilon}{2} |\oc_h(t)|_{H^1_*}^2 & + \int_0^t 2 \eta \|\obu_h\|_{\bH^1_*}^2 + (2 \eta)^{-1} s_0^{ext}(\op_h)^2 + \sigma M |\omu_h|_{H^1_*}^2 \ds \\
                    & = \int_0^t \sigma b(\obu_{h, \perp}, \oc_h, \wmu_h) - \sigma b(\wbu_{h, \perp}, \oc_h, \omu_h) - \sigma \varepsilon^{-1} (\ophi_h, \pa_t \oc_h)_* \ds \\
                    & \quad + \int_0^t \sigma \sD_\mu(\pa_t \oc_h) - \sD_\bu(\obu_h) - \sigma \sD_c(\omu_h) \ds.
                \end{split}
            \end{equation}
            To handle $(\ophi_h, \pa_t \oc_h)_*$ in \eqref{eq:err-eqn-strong}, we proceed as we did in \cite[Proof of Theorem~4.4]{OlshanskiiShakipov2026} by testing \eqref{eq:err-eq} with $v_h = \ophi_h \in V_h$:
            \[
                (\pa_t \oc_h, \ophi_h)_* = - M a_*(\omu_h, \ophi_h) - b(\obu_{h, \perp}, c_h, \ophi_h) - b(\wbu_{h, \perp}, \oc_h, \ophi_h) - \sD_c(\ophi_h).
            \]
            The other problematic term in \eqref{eq:err-eqn-strong} is $\sD_\mu(\pa_t \oc_h)$, which we handle using integration by parts in time:
            \[
                \int_0^t \sD_\mu(\pa_t \oc_h) \ds = \sD_\mu(t)(\oc_h(t)) - \int_0^t (\pa_t \sD_\mu)(\oc_h) \ds,
            \]
            where $\pa_t \sD_\mu$ is given in \eqref{eq:pat-sD-mu}, and we used $\oc_h(0) = 0$. Altogether,
            \begin{equation}
                \begin{split}
                    \frac{\sigma \varepsilon}{2} |\oc_h(t)|_{H^1_*}^2 & + \int_0^t 2 \eta \|\obu_h\|_{\bH^1_*}^2 + (2 \eta)^{-1} s_0^{ext}(\op_h)^2 + \sigma M |\omu_h|_{H^1_*}^2 \ds \\
                    & = \int_0^t \sigma b(\obu_{h, \perp}, \oc_h, \wmu_h) - \sigma b(\wbu_{h, \perp}, \oc_h, \omu_h) \ds \\
                    & \quad + \int_0^t [\sigma \varepsilon^{-1} M a_*(\omu_h, \ophi_h) + \sigma \varepsilon^{-1} b(\obu_{h, \perp}, c_h, \ophi_h) + \sigma \varepsilon^{-1} b(\wbu_{h, \perp}, \oc_h, \ophi_h)] \ds \\
                    & \quad + \sigma \sD_\mu(t)(\oc_h(t)) - \int_0^t [\sigma (\pa_t \sD_\mu)(\oc_h) + \sD_\bu(\obu_h) + \sigma \sD_c(\omu_h) - \sigma \varepsilon^{-1} \sD_c(\ophi_h)] \ds.
                \end{split}
            \end{equation}
            Testing \eqref{eq:err-eq} with $g_h = 1$ and using Lemma \ref{lem:ophih} yields
            \[
                |\mean{\omu_h}| \leq C_1 (|\oc_h|_{H^1_*}^2 + \|\sD_\mu\|_{V_h'}^2)^{1/2}.
            \]
            Thus, by \eqref{eq:discr-poincare}, we have
            \[
                \begin{split}
                    \frac{\sigma \varepsilon}{2} |\oc_h(t)|_{H^1_*}^2 & + \int_0^t 2 \eta \|\obu_h\|_{\bH^1_*}^2 + (2 \eta)^{-1} s_0^{ext}(\op_h)^2 + \sigma M \|\omu_h\|_{H^1_*}^2 \ds \\
                    & \leq C_2 \int_0^t |\oc_h|_{H^1_*}^2 + \|\sD_\mu\|_{V_h'}^2 \ds + \int_0^t \big[\underbrace{\sigma |b(\obu_{h, \perp}, \oc_h, \wmu_h)|}_{=: \rJ_1} + \underbrace{\sigma |b(\wbu_{h, \perp}, \oc_h, \omu_h)|}_{=: \rJ_2} \big] \ds \\
                    & + \int_0^t \big[ \underbrace{\sigma \varepsilon^{-1} M |a_*(\omu_h, \ophi_h)|}_{=: \rJ_3} + \underbrace{\sigma \varepsilon^{-1} |b(\obu_{h, \perp}, c_h, \ophi_h)|}_{=: \rJ_4} + \underbrace{\sigma \varepsilon^{-1} |b(\wbu_{h, \perp}, \oc_h, \ophi_h)|}_{=: \rJ_5} \big] \ds \\
                    & + \underbrace{\sigma |\sD_\mu(t)(\oc_h(t))|}_{=: \rJ_6} + \int_0^t \big[ \underbrace{\sigma |(\pa_t \sD_\mu)(\oc_h)|}_{=: \rJ_7} + \underbrace{|\sD_\bu(\obu_h)|}_{=: \rJ_8} + \underbrace{\sigma |\sD_c(\omu_h)|}_{=: \rJ_9} + \underbrace{\sigma \varepsilon^{-1} |\sD_c(\ophi_h)|}_{=: \rJ_{10}} \big] \ds.
                \end{split}
            \]
            We now estimate the ten terms. We first estimate the trilinear terms using the continuity from Lemma \ref{lem:b-disc}:
            \[
                \begin{split}
                    \rJ_1 & \leq C \|\obu_h\|_{\bH^1_*} |\oc_h|_{H^1_*} \|\wmu_h\|_{H^1_*} \leq \omega_1 \|\obu_h\|_{\bH^1_*}^2 + C_{\omega, 1} |\oc_h|_{H^1_*}^2 \|\wmu_h\|_{H^1_*}^2, \\
                    \rJ_2 & \leq C \|\wbu_h\|_{\bH^1_*} |\oc_h|_{H^1_*} \|\omu_h\|_{H^1_*} \leq \omega_2 \|\omu_h\|_{H^1_*}^2 + C_{\omega, 2} |\oc_h|_{H^1_*}^2 \|\wbu_h\|_{\bH^1_*}^2, \\
                    \rJ_4 & \leq C \|\obu_h\|_{\bH^1_*} \|c_h\|_{H^1_*} \|\ophi_h\|_{H^1_*} \leq \omega_4 \|\obu_h\|_{\bH^1_*}^2 + C_{\omega, 4} \|c_h\|_{L^\infty(I; H^1_*)}^2 \|\ophi_h\|_{H^1_*}^2, \\
                    \rJ_5 & \leq C \|\wbu_h\|_{\bH^1_*} |\oc_h|_{H^1_*} \|\ophi_h\|_{H^1_*} \leq C_{\rJ, 5} \|\wbu_h\|_{\bH^1_*}^2 |\oc_h|_{H^1_*}^2 + \|\ophi_h\|_{H^1_*}^2.
                \end{split}
            \]
            We next bound $\rJ_3$ and the defect terms $\rJ_6$-$\rJ_{10}$:
            \[
                \begin{aligned}
                    & \rJ_3 \leq \omega_3 \|\omu_h\|_{H^1_*}^2 + C_{\omega, 3} \|\ophi_h\|_{H^1_*}^2, \qquad & \rJ_6 & \leq \omega_6 |\oc_h(t)|_{H^1_*}^2 + C_{\omega, 6} \|\sD_\mu(t)\|_{V_h'}^2, \\
                    & \rJ_7 \leq C_{\rJ, 7} |\oc_h|_{H^1_*}^2 + \|\pa_t \sD_\mu\|_{V_h'}^2, \qquad & \rJ_8 & \leq \omega_8 \|\obu_h\|_{\bH^1_*}^2 + C_{\omega, 8} \|\sD_\bu\|_{\bV_h'}^2, \\
                    & \rJ_9 \leq \omega_9 \|\omu_h\|_{H^1_*}^2 + C_{\omega, 9} \|\sD_c\|_{V_h'}^2, \qquad & \rJ_{10} & \leq C_{\rJ, 10} \|\sD_c\|_{V_h'}^2 + \|\ophi_h\|_{H^1_*}^2.
                \end{aligned}
            \]
            Choosing $\omega_i > 0$ sufficiently small and using Lemma \ref{lem:ophih} to control $\ophi_h$ yields
            \[
                \begin{split}
                    |\oc_h(t)|_{H^1_*}^2+ \int_0^t \|\obu_h\|_{\bH^1_*}^2 + (2 \eta)^{-1} & s_0^{ext}(\op_h)^2 + \|\omu_h\|_{H^1_*}^2 \ds \leq C_1' \int_0^t \big( 1 + C_{\phi, 1}(s) + \|\wmu_h\|_{H^1_*}^2 + \|\wbu_h\|_{\bH^1_*}^2 \big) |\oc_h|_{H^1_*}^2 \ds \\
                    & \quad + C_3' \|\sD_\mu(t)\|_{V_h'}^2 + C_4' \int_0^t \big[ \|\sD_\mu\|_{V_h'}^2 + \|\sD_c\|_{V_h'}^2 + \|\sD_\bu\|_{\bV_h'}^2 + \|\pa_t \sD_\mu\|_{V_h'}^2 \big] \ds.
                \end{split}
            \]
            Since $C_{\phi, 1}(s) + \|\wmu_h\|_{H^1_*}^2 + \|\wbu_h\|_{\bH^1_*}^2 \in L^1(I)$, we deduce by Gr\"onwall's inequality that for a.e. $t \in I$
            \[
                \begin{split}
                    |\oc_h(t)|_{H^1_*}^2 + \int_0^t \|\obu_h\|_{\bH^1_*}^2 + (2 \eta)^{-1} s_0^{ext}(\op_h)^2 & + \|\omu_h\|_{H^1_*}^2 \ds \leq C_T^2 \left( \|\sD_\mu\|_{H^1(I; V_h')}^2 + \|\sD_c\|_{L^2(I; V_h')}^2 + \|\sD_\bu\|_{L^2(I; \bV_h')}^2 \right),
                \end{split}
            \]
            where we used $H^1(I; V_h') \hookrightarrow L^\infty(I; V_h')$ to simplify the $\sD_\mu$-term. This finishes the proof.
        \end{proof}
        Since $\oc_h\in V_{h,\#}$, discrete Poincar\'e \eqref{eq:discr-poincare} shows that \eqref{eq:disc-energy-est} also bounds $\|\oc_h\|_{L^\infty(I;H^1_*)}$ by $CC_T\cE_\sD$. From Theorem~\ref{thm:disc-energy} we deduce several further discrete error estimates collected in the corollary below.
        \begin{corollary} \label{cor:more-discr-err}
            Under the conditions of Theorem \ref{thm:disc-energy}, there holds
            \[
                \|\op_h\|_{L^2(I; L^2_*)} + \|\omu_h\|_{L^4(I; L^2_*)} \leq C C_T \cE_\sD.
            \]
            Furthermore, if $d=2$, then
            \[
                \|\oc_h\|_{L^\infty(I; L^\infty(\Omega_h))} \leq C |\log h|^{1/2} C_T \cE_\sD.
            \]
        \end{corollary}
        \begin{proof}
            We proceed in several steps: \smallskip

            \noindent \textbf{Pressure error}. Using the inf-sup condition \eqref{eq:discr-inf-sup} and the momentum error equation in \eqref{eq:err-eq}, we find
            \[
                \begin{split}
                    \alpha \|\op_h\|_{L^2_*} & \leq \sup_{\bzero \neq \bw_h \in \bV_h} \frac{(\nablaG \op_h, \bw_h)}{\|\bw_h\|_{\bH^1_*}} + \beta s_0^{ext}(\op_h) \\
                    & \leq 2 \eta \|\obu_h\|_{\bH^1_*} + \sigma C_b (\|c_h\|_{H^1_*} \|\omu_h\|_{H^1_*} + \|\oc_h\|_{H^1_*} \|\wmu_h\|_{H^1_*}) + \|\sD_\bu\|_{\bV_h'} + \beta s_0^{ext}(\op_h)
                \end{split}
            \]
           Squaring and integrating  both sides over $I$ yields
            \[
                \begin{split}
                    \|\op_h\|_{L^2(I; L^2_*)} & \leq C (\|\obu_h\|_{L^2(I; \bH^1_*)} + \|\omu_h\|_{L^2(I; H^1_*)} + \|\oc_h\|_{L^\infty(I; H^1_*)} + \|\sD_\bu\|_{L^2(I; \bV_h')} + \|s_0^{ext}(\op_h)\|_{L^2(I)}).
                \end{split}
            \]

            \noindent \textbf{Bound for $\omu_h$}. From the chemical potential error equation in \eqref{eq:err-eq} we have for a.e. $t \in I$,
            \[
                \|\omu_h(t)\|_{V_h'} \leq C (|\oc_h(t)|_{H^1_*} + \|\ophi_h(t)\|_{L^2_*} + \|\sD_\mu(t)\|_{V_h'}) \leq C (|\oc_h|_{L^\infty(I; H^1_*)} + \|\sD_\mu\|_{H^1(I; V_h')}).
            \]
            {Testing the discrete dual norm with $\omu_h$ gives $\|\omu_h\|_{L^2_*}^2\leq\|\omu_h\|_{V_h'}\|\omu_h\|_{H^1_*}$. Squaring and integrating yields the $L^4(I;L^2_*)$ estimate from the $L^\infty(I;V_h')$ and $L^2(I;H^1_*)$ bounds.} \smallskip
            
            \noindent \textbf{Pointwise bound for $d=2$}. This follows directly from Lemma \ref{lem:discr-sharp-Sobolev} and the error bound for $\|\oc_h\|_{L^\infty(I; H^1_*)}$. This finishes the proof.
        \end{proof}
        
        Theorem \ref{thm:disc-energy} bounds discrete errors by defects. The next lemma estimates the latter.
        \begin{lemma} \label{lem:defects}
            Let $(k_p+1, k_p, k_c) = (\deg \bV_h, \deg Q_h, \deg V_h)$, and let $\G$ be of class $C^{\max\{k_c + 2, k_p + 3\}}$. Assume further that
            \begin{equation} \label{eq:defects-reg}
                    c, \mu \in H^1(I; H^{k_c+1}(\G)), \qquad \bu \in L^\infty(I; \bH^2(\G)) \cap L^2(I; \bH^{k_p+2}(\G)), \qquad p \in L^2(I; H^{k_p+1}(\G)).
            \end{equation}
            Then
            \begin{equation} \label{eq:defect-bounds}
                \|\sD_\bu\|_{L^2(I; \bV_h')} + \|\sD_\mu\|_{H^1(I; V_h')} \leq C h^{k_c+1}, \qquad \|\sD_c\|_{L^2(I; V_h')} \leq C h^{\min\{k_c+1, k_p+2\}},
            \end{equation}
            where $C > 0$ depends only on $\G, \eta, \eta^{-1}, \sigma, M, L', \varepsilon^{-1}, \alpha$, shape-regularity of the triangulation, and the norms of $\bu, p, c, \mu$ bounded by the regularity assumption in \eqref{eq:defects-reg}.
        \end{lemma}
        \begin{proof}
            We estimate each defect separately. \smallskip
            
            \noindent \textit{Bound on $\sD_\bu$}. We estimate directly: since the $H^1$-projection preserves the mean, $\mu-\wmu_h$ and $c-\wc_h$ have zero mean and $\zm{\wmu_h}=R_h\zm{\mu}$. Let $\bw_{h, \perp} := \cP_E^\perp \bw_h \in \bH^1_{t, \perp}(\G)$, then
            \[
            \begin{split}
                \sigma^{-1} \sD_\bu(\bw_h) & = b(\bw_{h, \perp}, c, \mu) - b(\bw_{h, \perp}, \wc_h, \wmu_h) = b(\bw_{h, \perp}, c, \mu - \wmu_h) + b(\bw_{h, \perp}, c - \wc_h, \wmu_h) \\
                & = b(\bw_{h, \perp}, c, \mu - \wmu_h) - b(\bw_{h, \perp}, \zm{\wmu_h}, c - \wc_h) - (\divG \bw_{h, \perp}, \zm{\wmu_h} [c - \wc_h]) \\
                & = \underbrace{b(\bw_{h, \perp}, c, \mu - \wmu_h)}_{=: \fC_1}  \underbrace{- b(\bw_{h, \perp}, \zm{\mu}, c - \wc_h) - (\divG \bw_{h, \perp}, \zm{\mu} [c - \wc_h])}_{=: \fC_2} \\
                & \quad + \underbrace{b(\bw_{h, \perp}, \mu - \wmu_h, c - \wc_h) + (\divG \bw_{h, \perp}, [\mu - \wmu_h] [c - \wc_h])}_{=: \fC_3}.
            \end{split}
            \]
            Then
            \[
            \begin{split}
                |\fC_1| & \leq \|\bw_{h, \perp}\|_{\bH^1(\G)} \|(\nablaG c) (\mu - \wmu_h)\|_{\bH^1(\G)'} \leq \Csob C_K \|\bw_h\|_{\bH^1_*} \|\nablaG c\|_{L^3(\G)} \|\mu - \wmu_h\|, \\
                |\fC_2| & \leq C_K \|\bw_h\|_{\bH^1_*} \|(\nablaG \mu) (c - \wc_h)\|_{\bH^1(\G)'} + \|\divG \bw_{h, \perp}\| \|\zm{\mu}\|_{L^\infty(\G)} \|c - \wc_h\| \\
                & \leq 3 C_K \max\{\sqrt{d}, \Csob\} \|\bw_h\|_{\bH^1_*} \|\mu\|_{W^{1,3}(\G) \cap L^\infty(\G)} \|c - \wc_h\|, \\
                |\fC_3| & \leq C_K \|\bw_h\|_{\bH^1_*} \|[\nablaG (\mu - \wmu_h)] (c - \wc_h)\|_{\bH^1(\G)'} + \|\divG \bw_{h, \perp}\| \|\mu - \wmu_h\|_{L^4(\G)} \|c - \wc_h\|_{L^4(\G)} \\
                & \leq 2 C_K \max\{1, \sqrt{d} \Csob^2\} \|\bw_h\|_{\bH^1_*} \|\mu - \wmu_h\|_{H^1(\G)} \|c - \wc_h\|_{H^1(\G)}.
            \end{split}
            \]
            Above, we used Korn's inequality \eqref{eq:korn} to deduce $\|\nablaG \bw_{h,\perp}\| \leq \|\bw_{h, \perp}\|_{\bH^1(\G)} \leq C_K \|D_\G \bw_{h, \perp}\| \leq C_K \|\bw_h\|_{\bH^1_*}$, as well as $\|\divG \bv\| \leq \sqrt{d} \|\nablaG \bv\|$ for all $\bv \in \bH^1_t(\G)$. The intersection norm can be defined to be the maximum of the two norms. We thus arrive at
            \[
            \begin{split}
                \sigma^{-1} |\sD_\bu(\bw_h)| & \leq 3 C_K \max\{\sqrt{d}, \Csob\} \|\bw_h\|_{\bH^1_*} (\|\nablaG c\|_{L^3(\G)} + \|\mu\|_{W^{1,3}(\G) \cap L^\infty(\G)}) (\|\mu - \wmu_h\| + \|c - \wc_h\|) \\
                & \quad + 2 C_K \max\{1, \sqrt{d} \Csob^2\} \|\bw_h\|_{\bH^1_*} \|\mu - \wmu_h\|_{H^1(\G)} \|c - \wc_h\|_{H^1(\G)} \\
                & \leq C_3(t) \|\bw_h\|_{\bH^1_*} (\|\mu - \wmu_h\| + \|c - \wc_h\| + \|\mu - \wmu_h\|_{H^1(\G)} \|c - \wc_h\|_{H^1(\G)}),
            \end{split}
            \]
            where 
            \[
            C_3(t) := C' C_K (1 + \|\nablaG c\|_{L^3(\G)} + \|\mu\|_{W^{1,3}(\G) \cap L^\infty(\G)}), \quad C' = C'(\Csob, d).
            \]
            Note that $C_3(t) \in L^\infty(I)$ by assumption. We divide both sides above by $\|\bw_h\|_{\bH^1_*}$ and take supremum over $\bzero \neq \bw_h \in \bV_h$ to get
            \[
                \|\sD_\bu\|_{\bV_h'} \leq \sigma C_3(t) (\|\mu - \wmu_h\| + \|c - \wc_h\| + \|\mu - \wmu_h\|_{H^1(\G)} \|c - \wc_h\|_{H^1(\G)}).
            \]
            After squaring both sides and integrating in time, we deduce
            \[
            \begin{split}
                \int_I \|\sD_\bu\|_{\bV_h'}^2 \dt & \leq 3 \sigma^2 \|C_3\|_{L^\infty(I)}^2 \int_I \|\mu - \wmu_h\|^2 + \|c - \wc_h\|^2 + \|\mu - \wmu_h\|_{H^1(\G)}^2 \|c - \wc_h\|_{H^1(\G)}^2 \dt,
            \end{split}
            \]
            whence Lemma \ref{lem:prop-Rh} yields
            \[
            \begin{split}
                \|\sD_\bu\|_{L^2(I; \bV_h')} & \leq \sqrt3 \sigma \|C_3\|_{L^\infty(I)} \max\{\Cint, \Cint^2\} h^{k_c + 1} \\
                & \quad \times (\|\mu\|_{L^2(I; H^{k_c+1}(\G))} + \|c\|_{L^2(I; H^{k_c+1}(\G))} + \|c\|_{L^\infty(I; H^{k_c+1}(\G))} \|\mu\|_{L^2(I; H^{k_c+1}(\G))}),
            \end{split}
            \]
            where we used $h^{2k_c} \leq h^{k_c + 1}$ for $k_c \geq 1$. \smallskip

            \noindent \textit{Bound on $\sD_c$}. Since $c \in H^1(I; H^{k_c+1})$, we can write $(\pa_t \wt c_h, v_h)_* - \langle \pa_t c, v_h \rangle = (\pa_t [\wt c_h - c], v_h)_*$. Let $\wbu_{h,\perp} := \cP_E^\perp \wbu_h$. We have
            \[
            \begin{split}
                \sD_c(v_h) & = (\pa_t [\wc_h - c], v_h)_* + b(\wbu_{h,\perp}, \wc_h, v_h) - b(\bu, c, v_h) + M (\mu - \wmu_h, v_h)_* \\
                & = (\pa_t [\wc_h - c], v_h)_* + b(\wbu_{h,\perp} - \bu, \wc_h, v_h) + b(\bu, \wc_h - c, v_h) + M (\mu - \wmu_h, v_h)_* \\
                & = (\pa_t [\wc_h - c], v_h)_* + b(\wbu_{h,\perp} - \bu, c, v_h) + b(\wbu_{h,\perp} - \bu, \wc_h - c, v_h) + b(\bu, v_h, c - \wc_h) + M (\mu - \wmu_h, v_h)_*,
            \end{split}
            \]
            where at the end we used the fact that $\divG \bu = 0$ and so $b$ is skew-symmetric in its second and third arguments. Since the zero-mean operation in the third argument of $b$ is bounded in $H^1(\G)$, and $c - \wc_h$ has zero mean, whence $\zm{(c - \wc_h)} = c - \wc_h$, we estimate the terms as before:
            \[
            \begin{split}
                |\sD_c(v_h)| & \leq \|\pa_t [\wc_h - c]\|_{L^2_*} \|v_h\|_{L^2_*} + \Csob \|\wbu_{h, \perp} - \bu\| \|\nablaG c\|_{L^3(\G)} \|v_h\|_{H^1(\G)} \\
                & \quad + \Csob^2 \|\wbu_{h, \perp} - \bu\|_{\bH^1(\G)} \|\wc_h - c\|_{H^1(\G)} \|v_h\|_{H^1(\G)} + \|\bu\|_{L^\infty(\G)} \|v_h\|_{H^1(\G)} \|c - \wc_h\| + M \|\mu - \wmu_h\|_{L^2_*} \|v_h\|_{L^2_*} \\
                & \leq C_4(t) (\|\pa_t [\wc_h - c]\|_{L^2_*} + \|\wbu_h - \bu\|_{\bL^2_*} + \|\wbu_h - \bu\|_{\bH^1_*} \|\wc_h - c\|_{H^1(\G)} + \|c - \wc_h\| + \|\mu - \wmu_h\|_{L^2_*}) \|v_h\|_{H^1_*},
            \end{split}
            \]
            where $C_4(t) := \max \{1, M, \max\{1, \Csob\} (\|\nablaG c\|_{L^3(\G)} + \|\bu\|_{L^\infty(\G)}), \Csob^2 C_K\} \in L^\infty(I)$ by assumption (for $c$, we use the embedding $H^1(I) \hookrightarrow L^\infty(I)$ valid in one (time) dimension), and we bounded $H^1(\G)'$-norms by the $L^2(\G)$-norm for simplicity. Upon dividing both sides by $\|v_h\|_{H^1_*}$ and taking supremum over $v_h \neq 0$, we readily deduce from Lemmas \ref{lem:prop-Rh} and \ref{lem:prop-Stokes-proj} that
            \[
            \begin{split}
                \|\sD_c\|_{L^2(I; V_h')} & \leq \sqrt5 \|C_4\|_{L^\infty(I)} \max\{\Cint, \Cint^2\} (h^{k_c+1} + h^{k_p+2}) \\
                & \quad \times (\|c\|_{H^1(I; H^{k_c+1}(\G))} + \|\mu\|_{L^2(I; H^{k_c+1}(\G))} + \|\bu\|_{L^2(I; \bH^{k_p+2}(\G))} + \|p\|_{L^2(I; H^{k_p+1}(\G))} \\
                & \qquad \quad + \|c\|_{L^\infty(I; H^{k_c+1}(\G))} [\|\bu\|_{L^2(I; \bH^{k_p+2}(\G))} + \|p\|_{L^2(I; H^{k_p+1}(\G))}]),
            \end{split}
            \]
            where we used $h^{k_p+1+k_c} \leq h^{k_c + 1}$ for $k_p \geq 0$. \smallskip

            \noindent \textit{Bound on $\sD_\mu$}. Recall from \eqref{eq:sD-c-mu} that $\sD_\mu(g_h) = (\wmu_h - \mu, g_h)_* + \varepsilon^{-1}(f_0'(c) - f_0'(\wc_h), g_h) - \varepsilon (c - \wc_h, g_h)_*$. Throughout we use $\|g_h\|_{L^2_*} \leq \|g_h\|_{H^1_*}$ and $\varepsilon \leq 1$. By \eqref{eq:f0-consequences}, the Sobolev embedding $H^2(\G) \hookrightarrow L^\infty(\G)$ and Lemma \ref{lem:Linfty-Rh},
            \begin{equation} \label{eq:f0-diff-exact}
                \|f_0'(c) - f_0'(\wc_h)\| \leq L' \big( 1 + \|c\|_{L^\infty(\G)}^2 + \|\wc_h\|_{L^\infty(\G)}^2 \big) \|c - \wc_h\| \leq C \big(1 + \|c\|_{H^2(\G)}^2 \big) \|c - R_h c\| ,
            \end{equation}
            so that Lemma \ref{lem:prop-Rh} with $m = k_c+1$ gives, for a.e. $t \in I$,
            \[
                \|\sD_\mu\|_{V_h'} \leq C h^{k_c+1} \big( 1 + \|c\|^2_{H^2(\G)} \big) \big( \|\mu\|_{H^{k_c+1}(\G)} + \varepsilon^{-1} \|c\|_{H^{k_c+1}(\G)} \big),
            \]
            which is bounded uniformly in $t$ by \eqref{eq:defects-reg} and $H^1(I) \hookrightarrow L^\infty(I)$. We also recall \eqref{eq:pat-sD-mu}:
            \[
                (\pa_t \sD_\mu)(g_h) = (\pa_t[\wmu_h - \mu], g_h)_* + \varepsilon^{-1} \big( f_0''(c) \pa_t [c - \wc_h] + [f_0''(c) - f_0''(\wc_h)] \pa_t \wc_h, \, g_h \big) - \varepsilon (\pa_t [c - \wc_h], g_h)_* .
            \]
            By \eqref{eq:f0-consequences} and Lemma \ref{lem:Linfty-Rh},
            \[
                \begin{split}
                    \| f_0''(c) \pa_t[c - \wc_h] \| & \leq C \big( 1 + \|c\|_{H^2(\G)}^2 \big) \|\pa_t [c - R_h c]\|, \\
                    \big\| [f_0''(c) - f_0''(\wc_h)] \pa_t \wc_h \big\| & \leq C \big( 1 + \|c\|_{H^2(\G)} \big) \|c - R_h c\| \, \|R_h \pa_t c\|_{L^\infty(\G)} \leq C \big( 1 + \|c\|_{H^2(\G)} \big) \|c - R_h c\| \, \|\pa_t c\|_{H^2(\G)}.
                \end{split}
            \]
            Therefore, by Lemma \ref{lem:prop-Rh} applied to $c, \mu, \pa_t c$ and $\pa_t \mu$ with $m = k_c+1$,
            \[
                \|\pa_t \sD_\mu\|_{V_h'} \leq C h^{k_c+1} \big( 1 + \|c\|_{H^2(\G)}^2 \big) \big( \|\pa_t \mu\|_{H^{k_c+1}(\G)} + \varepsilon^{-1}\|\pa_t c\|_{H^{k_c+1}(\G)} + \varepsilon^{-1} \|c\|_{H^{k_c+1}(\G)} \|\pa_t c\|_{H^2(\G)} \big).
            \]
            Squaring and integrating in time, and using that $c \in C^0([0,T]; H^{k_c+1}(\G))$ and that $\pa_t c$ and $\pa_t \mu$ belong to $L^2(I; H^{k_c+1}(\G)) \subseteq L^2(I; H^2(\G))$, gives \eqref{eq:defect-bounds}. This finishes the proof.
        \end{proof}

    \subsection{Error estimates} \label{subsec:err-est}
        We now have all the ingredients needed to derive error estimates for our semidiscrete method.
        \begin{theorem} \label{thm:error-est}
            Under the assumptions of Theorem \ref{thm:disc-energy}  and Lemma \ref{lem:defects}, with $(\deg \bV_h, \deg Q_h, \deg V_h) = (k_p+1, k_p, k_c)$, there is $C > 0$ independent of $h$ and of how $\Gamma$ cuts through the mesh such that
            \begin{equation} \label{eq:error-est}
                \begin{split}
                    \|c - c_h\|_{L^\infty(I; H^1_*)} + \|\mu - \mu_h\|_{L^2(I; H^1_*)} & \leq C h^{\min\{k_c, k_p+2\}}, \\
                    \|\bu - \bu_h\|_{L^2(I; \bH^1_*)} + \|p - p_h\|_{L^2(I; L^2_*)} & \leq C h^{\min\{k_c+1, k_p+1\}}, \\
                    \|c - c_h\|_{L^\infty(I; L^2_*)} + \|\mu - \mu_h\|_{L^4(I; L^2_*)} + \|\bu - \bu_h\|_{L^2(I; \bL^2_*)} & \leq C h^{\min\{k_c+1, k_p+2\}}.
                \end{split}
            \end{equation}
        \end{theorem}
        \begin{proof}
            Lemma \ref{lem:defects} gives $\cE_\sD \leq C h^{\min\{k_c+1, k_p+2\}}$, so that Theorem \ref{thm:disc-energy} and Corollary \ref{cor:more-discr-err} yield
            \begin{equation} \label{eq:disc-energy-rates}
                \begin{split}
                    \|\obu_h\|_{L^2(I; \bH^1_*)} + \|\op_h\|_{L^2(I; L^2_*)} + \|\oc_h\|_{L^\infty(I; H^1_*)} + \|\omu_h\|_{L^2(I; H^1_*)} + \|\omu_h\|_{L^4(I; L^2_*)} \leq C h^{\min\{k_c+1, k_p+2\}}.
                \end{split}
            \end{equation}
            The first two rows in \eqref{eq:error-est} follow from the triangle inequality together with the projection error estimates $\|c - R_h c\|_{L^\infty(I; H^1_*)} + \|\mu - R_h \mu\|_{L^2(I; H^1_*)} \leq C h^{k_c}$ of Lemma \ref{lem:prop-Rh} (which require $c, \mu \in L^\infty(I; H^{k_c+1}(\G))$, resp. $L^2(I; H^{k_c+1}(\G))$) and $\|\bu - \bS_h \bu\|_{L^2(I; \bH^1_*)} + \|p - S_h p\|_{L^2(I; L^2_*)} \leq C h^{k_p+1}$ of Lemma \ref{lem:prop-Stokes-proj}.

            The $L^2$-errors in \eqref{eq:error-est} follow from the triangle inequality, the trivial bound $\|\cdot\|_{L^2_*} \leq C \|\cdot\|_{H^1_*}$ (resp. $\|\cdot\|_{\bL^2_*} \leq C \|\cdot\|_{\bH^1_*}$) and the superconvergent estimate \eqref{eq:disc-energy-rates}. 
        \end{proof}
        Since the Stokes--Cahn--Hilliard model features a diffuse interface, it is desirable to have an optimal-order (up to a logarithm) pointwise error bound for the phase field. For $d=2$, we have the following result.
        \begin{lemma}[pointwise error bound for $d=2$] \label{lem:pointwise-err-est}
            For $d=2$ and under the assumptions of Theorem \ref{thm:disc-energy}  and Lemma \ref{lem:defects}, with $(\deg \bV_h, \deg Q_h, \deg V_h) = (k_p+1, k_p, k_c)$, the following error estimate for $c$ is valid:
            \[
                \|c - c_h\|_{L^\infty(I; L^\infty(\Omega_h))} \leq C_1 \|c - R_h c\|_{L^\infty(I; L^\infty(\Omega_h))} + C_2 |\log h|^{1/2} h^{\min\{k_c+1, k_p+2\}}.
            \]
        \end{lemma}
        \begin{proof}
            Combining Lemma \ref{lem:defects} with the pointwise error estimate from Corollary \ref{cor:more-discr-err} yields
            \[
                \|\oc_h\|_{L^\infty(I; L^\infty(\Omega_h))} \leq C_2 |\log h|^{1/2} h^{\min\{k_c+1, k_p+2\}}.
            \]
            The final estimate then follows from the triangle inequality.
        \end{proof}
        \begin{remark}[regularity assumption]
            For $d=2$, the assumptions \eqref{eq:defects-reg} follow from Proposition~\ref{prop:well-posedness} (surface Stokes--Cahn--Hilliard) with $m=\max\{k_c+3,k_p-1\}$. In particular, $\mu\in H^1(I;H^{k_c+1}(\G))$ follows from $c\in H^1(I;H^m(\G))\cap C([0,T];H^{m+2}(\G))$ by differentiating \eqref{eq:sch-chem} in time, and $\bu\in L^\infty(I;\bH^2(\G))$ follows from Stokes regularity. Since the corresponding norms are bounded in terms of $\|c_0\|_{H^{m+2}(\G)}$~\cite[Theorem~4.7]{OlshanskiiShakipov2026}, so are the constants in Theorem~\ref{thm:error-est}.
        \end{remark}
        \begin{remark}[extensions of analysis] \label{rmk:extensions}
            The error analysis carried out above adapts straightforwardly to lifted parametric FEM \cite{DziukElliott2013} or to fitted FEM in flat (2D or 3D) domains under suitable boundary conditions. The latter corresponds to the FE discretizations of the standard model-H without inertia. Our treatment of $d=3$ becomes particularly relevant for the latter. In flat domains, we are not aware of a finite element error analysis for the quasi-static system; the closest results concern a Cahn--Hilliard--Darcy--Stokes system with a time derivative of the velocity \cite{DiegelFengWise2015} and a MAC finite difference scheme \cite{LiShen2020}. 
        \end{remark}
        \begin{remark}[dependence on $\varepsilon$]\label{rmk:robustness}
            Theorem~\ref{thm:error-est} holds for every fixed $\varepsilon\in(0,1]$, with constants depending on $\varepsilon$ through Gr\"onwall's inequality, as in related flat-domain analyses~\cite{DiegelFengWise2015,LiShen2020}. For the Allen--Cahn and Cahn--Hilliard equations, a spectral estimate in place of Gr\"onwall's inequality yields constants polynomial in $\varepsilon^{-1}$~\cite{FengProhl2003,FengProhl2004}. To the best of our knowledge, extending this approach to the classical model~H remains open, even in flat domains; we do not pursue it here.
        \end{remark}

\begin{remark}[choice of polynomial degrees]\label{rmk:balanced-degrees}
The best-approximation orders are $h^{k_c}$ for the phase variables in $H^1$, $h^{k_p+1}$ for the velocity in $\bH^1$ and pressure in $L^2$, and $h^{k_c+1}$ and $h^{k_p+2}$ for the phase variables and velocity, respectively, in spatial $L^2$ norms. The bounds in \eqref{eq:error-est} attain all these orders simultaneously when $k_c=k_p+1$. In particular, $k_c\geq k_p+1$ is sufficient for the estimate to give the optimal velocity $L^2$ order $h^{k_p+2}$, while $k_c\leq k_p+1$ avoids a loss in the phase-variable $L^2$ orders. Thus $k_c=k_p+1$ is the balanced choice. Results of numerical experiments support the predicted orders for  $k_c=k_p+1$.
\end{remark}
        
\section{Numerical experiment}
    For the numerical experiments, we use the polynomial double-well potential $f_0(s) = \frac14 (s^2-1)^2$. For time discretization, we apply an implicit midpoint rule with the secant potential; see \cite{OlshanskiiShakipov2026}. The method is formally second order in time. Numerical identification of the Killing fields at the discrete level is described in detail in \cite[\S6.1]{OlshanskiiShakipov2026}. The numerical method is implemented in \texttt{NGSolve/ngsxfem} \cite{Schoberl2009, Lehrenfeld2021}, which supports higher-order numerical approximation of geometry. 
    
    \subsection{Experimental order of convergence}
        The surface is given by a 2-1 prolate spheroid defined as the zero level set of $\phi(x, y, z) = x^2 + y^2 + z^2/4 - 1$ which has the one-dimensional Killing space $E = \bP \operatorname{span} \{(-y, x, 0)^T\}$. The PDE parameters are chosen to be $\eta = \sigma = M = \varepsilon = T = 1$. The manufactured solution is fixed to be 
        \[
        \begin{split}
            \wh \bu(x,y,z,t) & := (t + 1)\, \bn \times \nablaG [x(y + z^2)], \quad \bu := \cP_E^\perp \wh \bu, \quad \bu \cdot \bn = 0, \quad \divG \bu = 0, \\
            p(x,y,z,t) & := \cos (t) \sin(\pi x) \sin(\pi y z), \quad (p, 1) = 0, \\
            c(x,y,z,t) & := \frac{yz}{5} + \frac{\sin(t + x)}{2}, \\
            \mu & := \varepsilon^{-1} f_0'(c) - \varepsilon \DeltaG c.
        \end{split}
        \]
        This manufactured solution solves a forced version of \eqref{eq:weak}:
        \[
        \begin{split}
            - 2 \eta \Delta_S \bu + \nablaG p & = \sigma \zm{\mu} \nablaG c + \mathbf{f}, \\
            \divG \bu & = 0, \\
            \pa_t c + \bu \cdot \nablaG c - M \DeltaG \mu & = g, \\
            \mu & = \varepsilon^{-1} f_0'(c) - \varepsilon \DeltaG c,
        \end{split}
        \]
        where $\mathbf{f} := - 2 \eta \Delta_S \bu + \nablaG p - \sigma \zm{\mu} \nablaG c$ and $g := \pa_t c + \bu \cdot \nablaG c - M \DeltaG \mu$. 
        
        The $\bL^2(\G_h)$-part of the velocity error is computed as $\|\cP_h^\perp \bP_h (\bu^e - \bu_h)\|_{\G_h}$ where $\G_h$ is the approximate surface with piecewise smooth unit normal $\bn_h$, $\bP_h := \bI - \bn_h \otimes \bn_h$, and $\cP_h^\perp$ is the projection onto the orthogonal complement of discrete Killing fields as in \cite[\S5.1]{OlshanskiiShakipov2026}. For scalar-valued quantities we simply compute $(\|v^e - v_h\|_{\G_h}^2 + h \|\bn_h \cdot \nabla v_h\|_{\Omega_h}^2)^{1/2}$. The results along with further details are given in Table \ref{tab:conv}.

        \begin{table}[ht]
            \centering
            \begin{tabular}{c|c|cc|cc|cc|cc|cc}
            \hline
            \multirow{2}{*}{$h$} & \multirow{2}{*}{dof} & \multicolumn{2}{c|}{$\|e_{\mathbf{u}}\|_{L^2(I; \bL^2_*)}$} & \multicolumn{2}{c|}{$\|e_{\mathbf{u}}\|_{L^2(I; \bH^1_*)}$} & \multicolumn{2}{c|}{$\|e_p\|_{L^2(I; L^2_*)}$} & \multicolumn{2}{c|}{$\|e_c\|_{L^\infty(I; L^\infty(\G_h))}$} & \multicolumn{2}{c|}{$\|e_\mu\|_{L^4(I; L^2_*)}$} \\
            \cline{3-12}
             & & error & rate & error & rate & error & rate & error & rate & error & rate \\
            \hline
            \multicolumn{12}{c}{$(\bP_2, P_1, P_1, P_1)$, $\Delta t = h^{3/2}$} \\
            \hline
            $\mathtt{3.13e\text{-}1}$ & $\mathtt{12 257}$ & $\mathtt{7.04e\text{-}1}$ & -- & $\mathtt{2.33e0}$ & -- & $\mathtt{8.47e\text{-}1}$ & -- & $\mathtt{1.52e\text{-}1}$ & -- & $\mathtt{5.93e\text{-}1}$ & -- \\
            $\mathtt{1.56e\text{-}1}$ & $\mathtt{49 559}$ & $\mathtt{6.57e\text{-}2}$ & $\mathtt{3.42}$ & $\mathtt{5.12e\text{-}1}$ & $\mathtt{2.19}$ & $\mathtt{2.13e\text{-}1}$ & $\mathtt{1.99}$ & $\mathtt{2.88e\text{-}2}$ & $\mathtt{2.40}$ & $\mathtt{1.45e\text{-}1}$ & $\mathtt{2.03}$ \\
            $\mathtt{7.81e\text{-}2}$ & $\mathtt{196 637}$ & $\mathtt{8.51e\text{-}3}$ & $\mathtt{2.95}$ & $\mathtt{1.27e\text{-}1}$ & $\mathtt{2.01}$ & $\mathtt{5.08e\text{-}2}$ & $\mathtt{2.07}$ & $\mathtt{8.59e\text{-}3}$ & $\mathtt{1.75}$ & $\mathtt{3.95e\text{-}2}$ & $\mathtt{1.87}$ \\
            $\mathtt{3.91e\text{-}2}$ & $\mathtt{791 096}$ & $\mathtt{1.25e\text{-}3}$ & $\mathtt{2.77}$ & $\mathtt{3.13e\text{-}2}$ & $\mathtt{2.02}$ & $\mathtt{1.22e\text{-}2}$ & $\mathtt{2.05}$ & $\mathtt{2.33e\text{-}3}$ & $\mathtt{1.88}$ & $\mathtt{1.03e\text{-}2}$ & $\mathtt{1.94}$ \\
            \hline
            \multicolumn{12}{c}{$(\bP_2, P_1, P_2, P_2)$, $\Delta t = h^{3/2}$} \\
            \hline
            $\mathtt{3.13e\text{-}1}$ & $\mathtt{18 081}$ & $\mathtt{6.91e\text{-}1}$ & -- & $\mathtt{2.32e0}$ & -- & $\mathtt{8.36e\text{-}1}$ & -- & $\mathtt{1.51e\text{-}2}$ & -- & $\mathtt{1.66e\text{-}1}$ & -- \\
            $\mathtt{1.56e\text{-}1}$ & $\mathtt{73 027}$ & $\mathtt{6.60e\text{-}2}$ & $\mathtt{3.39}$ & $\mathtt{5.12e\text{-}1}$ & $\mathtt{2.18}$ & $\mathtt{2.12e\text{-}1}$ & $\mathtt{1.98}$ & $\mathtt{1.77e\text{-}3}$ & $\mathtt{3.09}$ & $\mathtt{2.04e\text{-}2}$ & $\mathtt{3.02}$ \\
            $\mathtt{7.81e\text{-}2}$ & $\mathtt{289 897}$ & $\mathtt{8.11e\text{-}3}$ & $\mathtt{3.02}$ & $\mathtt{1.27e\text{-}1}$ & $\mathtt{2.01}$ & $\mathtt{5.05e\text{-}2}$ & $\mathtt{2.07}$ & $\mathtt{2.81e\text{-}4}$ & $\mathtt{2.66}$ & $\mathtt{2.61e\text{-}3}$ & $\mathtt{2.96}$ \\
            $\mathtt{3.91e\text{-}2}$ & $\mathtt{1 166 230}$ & $\mathtt{1.01e\text{-}3}$ & $\mathtt{3.01}$ & $\mathtt{3.13e\text{-}2}$ & $\mathtt{2.02}$ & $\mathtt{1.22e\text{-}2}$ & $\mathtt{2.05}$ & $\mathtt{3.74e\text{-}5}$ & $\mathtt{2.91}$ & $\mathtt{3.24e\text{-}4}$ & $\mathtt{3.01}$ \\
            \hline
            \multicolumn{12}{c}{$(\bP_3, P_2, P_3, P_3)$, $\Delta t = h^{2}$} \\
            \hline
            $\mathtt{3.13e\text{-}1}$ & $\mathtt{55 938}$ & $\mathtt{2.31e\text{-}1}$ & -- & $\mathtt{9.27e\text{-}1}$ & -- & $\mathtt{3.13e\text{-}1}$ & -- & $\mathtt{8.18e\text{-}3}$ & -- & $\mathtt{5.20e\text{-}2}$ & -- \\
            $\mathtt{1.56e\text{-}1}$ & $\mathtt{224 996}$ & $\mathtt{2.57e\text{-}2}$ & $\mathtt{3.17}$ & $\mathtt{1.50e\text{-}1}$ & $\mathtt{2.63}$ & $\mathtt{4.33e\text{-}2}$ & $\mathtt{2.85}$ & $\mathtt{3.88e\text{-}4}$ & $\mathtt{4.40}$ & $\mathtt{4.14e\text{-}3}$ & $\mathtt{3.65}$ \\
            $\mathtt{7.81e\text{-}2}$ & $\mathtt{894 262}$ & $\mathtt{1.45e\text{-}3}$ & $\mathtt{4.14}$ & $\mathtt{1.68e\text{-}2}$ & $\mathtt{3.15}$ & $\mathtt{4.73e\text{-}3}$ & $\mathtt{3.19}$ & $\mathtt{2.43e\text{-}5}$ & $\mathtt{4.00}$ & $\mathtt{2.79e\text{-}4}$ & $\mathtt{3.89}$ \\
            \hline
            \end{tabular}
            \caption{Experimental order of convergence. We use $\Delta t = h^{(k_p+2)/2}$ (needed for $L^2$-error equilibration) on the 2-1 prolate spheroid. Here, $e_{\xi} := \xi - \xi_h$ are errors for $\xi \in \{\bu, p, c, \mu\}$ measured in different norms. The geometry approximation order used in the experiment is $\max\{k_p + 1, k_c\}$, which is consistent with \cite[Rem.~5.1]{OlshanskiiShakipov2026}. The starred norms are direct discrete analogs of the starred norms used in the analysis, i.e., they contain all stabilization terms.}
            \label{tab:conv}
        \end{table}
        
The results for the balanced (cf. Remark~\ref{rmk:balanced-degrees}) pairs $(\bP_2,P_1,P_2,P_2)$ and $(\bP_3,P_2,P_3,P_3)$ show velocity $L^2$ rates close to three and four, respectively, and velocity $H^1$ and pressure $L^2$ rates close to two and three. The chemical potential also converges with spatial $L^2$ orders close to three and four. These observations agree with Theorem~\ref{thm:error-est}. The phase-field maximum-norm errors exhibit rates close to $k_c+1$, consistent with Lemma~\ref{lem:pointwise-err-est} if the Ritz-projection error has the corresponding pointwise order. For $(\bP_2,P_1,P_1,P_1)$, the phase and chemical potential errors converge at approximately second order, whereas the last two velocity $L^2$ rates are $2.95$ and $2.77$. Note that the theorem established only $O(h^2)$ convergence rate for the case of  $k_c=k_p=1$. Further refinement would be needed to determine the experimental asymptotic rate in this case. 

\begin{acknowledgement}
AI tools (Claude Opus 5.5 and ChatGPT 6 Astra) assisted in extending the trace and normal-extension estimates to the $L^p$ setting and general dimension, proving the equivalence of the discrete inf-sup conditions, identifying the elliptic-reconstruction approach of Heywood and Rannacher~\cite{HeywoodRannacher1982} for the discrete Sobolev embedding, and developing numerical code. The authors independently verified the mathematical arguments, checked the numerical implementation and results, and take full responsibility for the content of the paper.
\end{acknowledgement}

\begin{funding}
    The authors were partially supported by {the} U.S. National Science Foundation under award DMS-2408978.
\end{funding}

\appendix

\section{Trace inequality and stability of the constant normal extension in the \texorpdfstring{$L^p$}{Lp} setting}
    A local trace inequality for a simplex $T$ cut by $\G$ is a fundamental result in the analysis of TraceFEM that gives estimates uniform with respect to the position of the cut. The $L^2$-based version was proven in \cite[Lemma 4.3]{Reusken2015} for surfaces. We note that this result does not follow from the classical trace inequality or divergence theorem because $T$ cut by a curved $\G$ can in principle yield a non-Lipschitz domain.

    Below we present the result in the $L^p$-based setting valid for any $d \geq 1$ and $p \in [1,\infty]$. The main idea is to notice that when the normal does not oscillate much on $T$, the cut $\G \cap T$ can be seen as a graph of a $C^1$ function over a hyperplane. The estimate follows by integrating along chords of $T$ perpendicular to this hyperplane and applying the standard trace inequality on the boundary of the convex polytope $T$.
    \begin{lemma}[local trace inequality in $L^p$] \label{lem:local-trace-ineq}
        Let $\G$ be a closed $d$-dimensional hypersurface of class $C^{1,1}$, $d \in \bbZ_{\geq 1}$, let $T \subseteq U_\delta$ be a convex polytope with diameter $h_T$ and inradius $\rho_T$, and assume that $h_T L_T \leq \frac12$, where $L_T := \|\nabla \bn\|_{\dL^\infty(T)}$. Then for all $p \in [1, \infty]$ and $v \in W^{1,p}(T)$,
        \[
            \|v\|_{L^p(\G_T)} \leq C \big(h_T^{-1/p} \|v\|_{L^p(T)} + h_T^{1-1/p} \|\nabla v\|_{\bL^p(T)}\big), \qquad \G_T := \G \cap T,
        \]
        where $C > 0$ depends only on $p$ and $h_T/\rho_T$, and in particular not on how $\G$ cuts $T$. Here $1/\infty := 0$.
    \end{lemma}
    \begin{proof}
        If $\G_T = \emptyset$ there is nothing to prove. For $p = \infty$, $v$ has a representative in $C^{0,1}(\overline T)$ because $T$ is convex, so $\|v\|_{L^\infty(\G_T)} \leq \|v\|_{L^\infty(T)}$. Let $p \in [1, \infty)$. Since $C^1(\overline{T})$ is dense in $W^{1,p}(T)$, it suffices to prove the estimate for $v \in C^1(\overline{T})$.

        \emph{Step 1: $\G_T$ is a graph.} All quantities in the lemma are invariant under rotations of $\bbR^{d+1}$, so we may fix $\bx_T \in \G_T$ and assume that $\bn(\bx_T) = \be_{d+1}$ (the last coordinate vector). We write points of $\bbR^{d+1}$ as $(\bz, s) \in \bbR^d \times \bbR$ and $\bn = (\bn', n_{d+1})$; see Figure~\ref{fig:local-trace}. Since $T$ is convex and $\bn$ is Lipschitz on $T$ with constant $L_T$, we have $|\bn - \be_{d+1}| \leq L_T h_T \leq \frac12$ on $T$, and since both are unit vectors,
        \begin{equation} \label{eq:monotone}
            n_{d+1} = 1 - \tfrac12 |\bn - \be_{d+1}|^2 \geq \tfrac78 \qquad \text{on } T.
        \end{equation}
        Let $\cV \subseteq \bbR^d$ be the projection of $T$ onto the first $d$ coordinates. For $\bz \in \cV$, the vertical line through $\bz$ meets the convex set $T$ in a segment $\{\bz\} \times [a(\bz), b(\bz)]$, which we call a chord. Along a chord the signed distance is strictly increasing, since $\pa_s \dst = n_{d+1} \geq \frac78$ by \eqref{eq:monotone}. Hence every chord meets $\G_T$, the zero set of $\dst$ in $T$, at most once, and $\G_T$ is the graph of a function $\tau$ over its projection $\cV_\G \subseteq \cV$:
        \[
            \G_T = \{(\bz, \tau(\bz)) : \bz \in \cV_\G\}.
        \]
        Since $\dst$ is of class $C^1$ and $\pa_s \dst > 0$ on $T$, the implicit function theorem shows that near each point of $\G_T$, the surface $\G$ is the graph $s = \psi(\bz)$ of a $C^1$ function with $\nabla \psi = -\bn'/n_{d+1}$. Hence, near each point of $\G_T$, the surface element is
        \[
            \dA = \sqrt{1 + |\nabla \psi|^2} \dz = \frac{\mathrm{d}\bz}{n_{d+1}}, \qquad \text{because} \qquad 1 + \frac{|\bn'|^2}{n_{d+1}^2} = \frac{|\bn|^2}{n_{d+1}^2} = \frac{1}{n_{d+1}^2}.
        \]
        Since this identity is local and the projection $(\bz, s) \mapsto \bz$ is injective on $\G_T$, we obtain for every nonnegative $w \in C^0(\overline{T})$
        \begin{equation} \label{eq:area-formula}
            \tfrac78 \int_{\G_T} w \dA \leq \int_{\G_T} w \, n_{d+1} \dA = \int_{\cV_\G} w(\bz, \tau(\bz)) \dz.
        \end{equation}

        \begin{figure}[htbp]
            \centering
            \begin{tikzpicture}[
                >=Stealth, line join=round, line cap=round, font=\small,
                body/.style   = {fill=black!5, draw=black!55},
                entry/.style  = {draw=red!60!black, line width=1.2pt},
                surf/.style   = {draw=green!40!black, line width=1.2pt},
                chord/.style  = {draw=black!75},
                chordup/.style= {draw=black!35, densely dotted},
                rayout/.style = {draw=black!45, dash pattern=on 2pt off 1.6pt},
                hplane/.style = {draw=black!55},
                Vseg/.style   = {draw=blue!55!black, line width=1pt},
                Wseg/.style   = {draw=orange!85!black, line width=2pt},
                cone/.style   = {fill=blue!8, draw=blue!30!black, line width=.3pt},
                adot/.style   = {circle, fill=red!60!black, inner sep=0pt, minimum size=3.2pt},
                tdot/.style   = {circle, fill=blue!55!black, inner sep=0pt, minimum size=3.4pt},
                leader/.style = {draw=black!45, line width=.35pt}
            ]
            \def\eang{90}
            \pgfmathsetmacro{\pang}{\eang-90}
            \def\zs#1#2#3{%
                \pgfmathsetmacro{\tmpx}{(#1)*cos(\pang)+(#2)*cos(\eang)}%
                \pgfmathsetmacro{\tmpy}{(#1)*sin(\pang)+(#2)*sin(\eang)}%
                \coordinate (#3) at (\tmpx,\tmpy);}
            \def\tauexpr{(2.15-0.11*(\zz-0.15)*(\zz-0.15))}
     
            \zs{-2.6}{0.7}{A} \zs{2.8}{1.0}{B} \zs{-0.1}{4.6}{C}
     
            \zs{-3.6}{0}{H1} \zs{3.8}{0}{H2}
            \zs{-2.6}{0}{V1} \zs{2.8}{0}{V2}
            \zs{-1.995}{0}{W1} \zs{2.273}{0}{W2}
     
            \zs{-1.6}{0}{bs1} \zs{-1.6}{0.756}{en1} \zs{-1.6}{1.813}{su1} \zs{-1.6}{2.260}{ex1}
            \zs{-0.6}{0}{bs2} \zs{-0.6}{0.811}{en2} \zs{-0.6}{2.088}{su2} \zs{-0.6}{3.820}{ex2}
            \zs{1.0}{0}{bs3}  \zs{1.0}{0.900}{en3}  \zs{1.0}{2.071}{su3}  \zs{1.0}{3.235}{ex3}
            \zs{1.8}{0}{bs4}  \zs{1.8}{0.944}{en4}  \zs{1.8}{1.851}{su4}  \zs{1.8}{2.241}{ex4}
     
            \zs{0.15}{2.15}{xT} \zs{0.15}{3.05}{eT}
            \zs{-0.2857}{2.9375}{k1} \zs{-0.0754}{3.0213}{k2}
            \zs{0.3754}{3.0213}{k3} \zs{0.5857}{2.9375}{k4}
     
            \zs{-1.2}{1.9495}{n0a} \zs{-1.3708}{2.5247}{n1a}
            \zs{1.35}{1.9916}{n0b} \zs{1.5031}{2.5717}{n1b}
     
            \zs{-2.05}{1.15}{labT} \zs{-2.9}{1.90}{labG} \zs{0.15}{1.88}{labX}
            \zs{0.15}{3.25}{labE} \zs{2.0}{3.35}{labC} \zs{1.5}{3.30}{labCa}
            \zs{0.2}{0}{labW}     \zs{2.56}{0}{labV}
            \zs{-1.995}{1.644}{gL}
     
            \draw[hplane] (H1) -- (H2) node[right] {$\bbR^d$};
            \draw[Vseg] (V1) -- (V2);
            \draw[Wseg] (W1) -- (W2);
     
            \filldraw[body] (A) -- (B) -- (C) -- cycle;
            \draw[entry] (A) -- (B);
     
            \foreach \i in {1,2,3,4}{%
                \draw[rayout]  (bs\i) -- (en\i);
                \draw[chord]   (en\i) -- (su\i);
                \draw[chordup] (su\i) -- (ex\i);
                \fill[black!55] (bs\i) circle (1pt);
            }
     
            \draw[cone] (xT) -- (k1) -- (k2) -- (eT) -- (k3) -- (k4) -- cycle;
            \draw[surf] plot[domain=-1.995:2.273,samples=80,variable=\zz,smooth]
                ({(\zz)*cos(\pang)+\tauexpr*cos(\eang)},
                 {(\zz)*sin(\pang)+\tauexpr*sin(\eang)});
            \draw[->] (n0a) -- (n1a);
            \draw[->] (n0b) -- (n1b);
            \draw[->,thick] (xT) -- (eT);
     
            \foreach \i in {1,2,3,4}{%
                \node[adot] at (en\i) {};
                \node[tdot] at (su\i) {};
            }
            \node[tdot] at (xT) {};
     
            \node at (labT) {$T$};
            \node (gam) at (labG) {$\G_T$};
            \draw[leader] (gam) -- (gL);
     
            \node[below] at (labX) {$\bx_T$};
            \node at (labE) {$\be_{d+1}$};
            \node at (labC) {$n_{d+1} \geq \tfrac78$};
            \draw[leader] (k4) -- (labCa);
            \node[below] at (labW) {$\cV_\G$};
            \node[below] at (labV) {$\cV$};
            \end{tikzpicture}
            \caption{Proof of Lemma~\ref{lem:local-trace-ineq}. Since $h_T L_T \leq \frac12$,
                the normal stays in a cone of half-angle $\arccos\frac78 \approx 29^\circ$ about
                $\be = \bn(\bx_T)$, so every line parallel to $\be$ meets $\G_T$ at most once and
                $\G_T$ is a graph over $\cV_\G = \pi(\G_T)$. Each point of $\G_T$ (blue) is traded
                for the entry point of its chord on $\pa T$ (red) plus a gradient integral along
                the chord; nothing here depends on where $\G$ cuts $T$.}
            \label{fig:local-trace}
        \end{figure}
        
        \emph{Step 2: slicing.} Let $w := |v|^p$ and $g := p|v|^{p-1}|\nabla v|$ (with $|v|^0 := 1$), both nonnegative and continuous on $\overline T$. Along every chord, $s \mapsto w(\bz, s)$ is Lipschitz with $|\pa_s w| \leq g$ a.e., and $|\nabla w| \leq g$ a.e. in $T$. Hence, for every $\bz \in \cV_\G$,
        \[
            w(\bz, \tau(\bz)) \leq w(\bz, a(\bz)) + \int_{a(\bz)}^{b(\bz)} g(\bz, s) \ds.
        \]
        We integrate over $\cV_\G$ and enlarge the domain to $\cV$, which is allowed since $w, g \geq 0$. By Fubini, the chord integrals add up to $\|g\|_{L^1(T)}$. The lower ends $(\bz, a(\bz))$ of the chords form the bottom of $\pa T$, the graph of the convex, piecewise affine function $a$. Its surface element $\sqrt{1 + |\nabla a|^2} \dz$ is at least $\mathrm{d}\bz$, so $\int_\cV w(\bz, a(\bz)) \dz \leq \|w\|_{L^1(\pa T)}$. We combine this with \eqref{eq:area-formula}, the classical trace inequality $\|w\|_{L^1(\pa T)} \leq C (h_T^{-1}\|w\|_{L^1(T)} + \|\nabla w\|_{\bL^1(T)})$, with $C$ depending only on $h_T/\rho_T$, and H\"older's inequality $\|g\|_{L^1(T)} \leq p \|v\|_{L^p(T)}^{p-1} \|\nabla v\|_{\bL^p(T)}$. This yields
        \begin{equation} \label{eq:trace-L1}
            \|v\|_{L^p(\G_T)}^p \leq C \big( h_T^{-1} \|v\|_{L^p(T)}^p + p \|v\|_{L^p(T)}^{p-1} \|\nabla v\|_{\bL^p(T)} \big).
        \end{equation}
        \emph{Step 3: scaling.} By Young's inequality with exponents $\frac{p}{p-1}$ and $p$ (trivial if $p = 1$),
        \[
            p \|v\|_{L^p(T)}^{p-1} \|\nabla v\|_{\bL^p(T)} \leq (p-1) h_T^{-1} \|v\|_{L^p(T)}^p + h_T^{p-1} \|\nabla v\|_{\bL^p(T)}^p.
        \]
        Inserting this into \eqref{eq:trace-L1} and taking the $p$-th root finishes the proof.
    \end{proof}
        
    The following corollary is an immediate consequence of Lemma \ref{lem:local-trace-ineq} and the local inverse inequality.
    \begin{corollary}[FE trace inequality] \label{cor:local-trace-discrete}
        Under the conditions of Lemma \ref{lem:local-trace-ineq}, any polynomial $v_h$  satisfies
        \[
            \|v_h\|_{L^p(\G_T)} \leq C h_T^{-1/p} \|v_h\|_{L^p(T)},
        \]
        where $C$ depends also on the polynomial degree of $v_h$.
    \end{corollary}

    The next result extends \cite[Lemma 3.1]{Reusken2015} to the $W^{k,p}$-setting, $p \in [1, \infty]$, and to any intrinsic dimension $d \geq 1$. The proof follows the same idea, with the case $p = \infty$ requiring a little more care. We also propose a simplified argument for $k \geq 2$. 
    \begin{lemma}[$W^{k,p}$-stability of the constant normal extension] \label{lem:const-normal-ext}
        Let $\G$ be a closed $d$-dimensional hypersurface of class $C^{\max\{2, k+1\}}$, $d \in \bbZ_{\geq 1}$, $k \in \bbZ_{\geq 0}$, let $\delta_0, \delta$ be as in Section~\ref{sec:semidiscrete} (the construction there works for any $d \geq 1$), and let $p \in [1, \infty]$. Then for all $r \in (0, \delta]$ and $v \in W^{k,p}(\G)$, we have $v^e \in W^{k,p}(U_r)$ and
        \[
            \|v^e\|_{W^{k,p}(U_r)} \leq C\, r^{1/p} \|v\|_{W^{k,p}(\G)},
        \]
        where $C > 0$ depends only on $d$ if $k \leq 1$, and on $d$, $k$ and $\|D^j \dst\|_{L^\infty(U_{\delta_0})}$, $j \leq k+1$, if $k \geq 2$.
    \end{lemma}
 
    \begin{proof}
        \emph{Normal coordinates.} The map $(\by, s) \mapsto \by + s \bn(\by)$ is a bijection $\G \times (-r, r) \to U_r$, under which $\dst(\by + s \bn(\by)) = s$ and $\bp(\by + s \bn(\by)) = \by$, and $\dx = J(\by, s) \ds \dA(\by)$ with $J(\by, s) = \prod_{i=1}^d (1 + s \kappa_i(\by))$. Moreover, the eigenvalues of $\nabla \nabla \dst$ at $\by + s\bn(\by)$ are $0$ and $\frac{\kappa_i(\by)}{1 + s \kappa_i(\by)}$, for $i = 1, \dots, d$ \cite[Lemmas 6 and 22]{BonitoDemlowNochetto2020}. Recall from \eqref{eq:assump-delta} that $5 \delta \kappa_{\max} < 1$, whence $|s| \kappa_{\max} \leq \delta \kappa_{\max} < \frac15$, and we obtain in $U_r$
        \begin{equation} \label{eq:jacobian-bounds}
            \big(\tfrac45\big)^d \leq J \leq \big(\tfrac65\big)^d, \qquad |\nabla \bn| = |\nabla \nabla \dst| \leq \tfrac54 \kappa_{\max}, \qquad |\bM| \leq \tfrac54,
        \end{equation}
        where $\bM := \bI - \dst \, \nabla \nabla \dst$ has eigenvalues $1$ and $\frac{1}{1 + s\kappa_i}$. In particular, for every measurable $u : \G \to \bbR$ and $p \in [1, \infty)$,
        \begin{equation} \label{eq:cne-change}
            \|u^e\|_{L^p(U_r)}^p = \int_{-r}^r \int_\G |u(\by)|^p J(\by, s) \dA(\by) \ds \leq 2 \big(\tfrac65\big)^d r \|u\|_{L^p(\G)}^p.
        \end{equation}
 
        \emph{Case $k \leq 1$.} For $k = 0$ and $p < \infty$, this is \eqref{eq:cne-change} with $u = v$. For $p = \infty$, the identity in \eqref{eq:cne-change} with $u = \mathbf{1}_N$ shows that $|\bp^{-1}(N)| = 0$ whenever $\mathcal H^d(N) = 0$; taking $N = \{|v| > \|v\|_{L^\infty(\G)}\}$ gives $\|v^e\|_{L^\infty(U_r)} \leq \|v\|_{L^\infty(\G)}$.
 
        For $k = 1$, recall $\nabla v^e = \bM (\nablaG v)^e$ in $U_r$ for $v \in C^1(\G)$ \cite[Eq.~(27)]{BonitoDemlowNochetto2020}. Let $p < \infty$, and take $v_m \in C^1(\G)$ with $v_m \to v$ in $W^{1,p}(\G)$ \cite[Theorem 3.6]{BenavidesNochettoShakipov2026-scalar}. By \eqref{eq:cne-change} and $|\bM| \leq \frac54$, $v_m^e \to v^e$ and $\nabla v_m^e \to \bM (\nablaG v)^e$ in $L^p(U_r)$, so $v^e \in W^{1,p}(U_r)$ and $\nabla v^e = \bM (\nablaG v)^e$. The bound follows from $|\nabla v^e| \leq \frac54 |\nablaG v|^e$ and \eqref{eq:cne-change} with $u = |\nablaG v|$. If $p = \infty$, the identity holds a.e. because $W^{1,\infty}(\G) \subseteq W^{1,q}(\G)$ for $q < \infty$, and the null-set argument above gives $\|\nabla v^e\|_{\bL^\infty(U_r)} \leq \frac54 \|\nablaG v\|_{\bL^\infty(\G)}$.
 
        \emph{Case $k \geq 2$ (sketch).} Write $\uD_l v := (\nablaG v)_l$, so that $\pa_j v^e = \bM_{jl} (\uD_l v)^e$ (summation over repeated indices). For $v \in C^k(\G)$ and $\G$ of class $C^{k+1}$, applying this identity once more to $\uD_l v$ gives
        \[
            \pa_i \pa_j v^e = (\pa_i \bM_{jl}) (\uD_l v)^e + \bM_{jl} \bM_{is} (\uD_s \uD_l v)^e,
        \]
        where $\pa_i \bM_{jl}$ involves third derivatives of $\dst$. Iterating gives $|D^j v^e| \leq C \sum_{i=1}^j |\nablaG^i v|^e$ in $U_r$ for $1 \leq j \leq k$, with $C$ depending on $d$, $k$ and $\|D^m \dst\|_{L^\infty(U_{\delta_0})}$, $m \leq k+1$. The claim follows from \eqref{eq:cne-change} with $u = |\nablaG^i v|$ and $\|\nablaG^i v\|_{L^p(\G)} \leq C \|v\|_{W^{k,p}(\G)}$ \cite[Proposition 3.7]{BenavidesNochettoShakipov2026-scalar}, and for general $v$ as in the case $k = 1$.
    \end{proof}

\bibliographystyle{dinat}
\bibliography{references}

@article {OlshanskiiReuskenZhiliakov2021,
    AUTHOR = {Olshanskii, Maxim A. and Reusken, Arnold and Zhiliakov,
              Alexander},
     TITLE = {Inf-sup stability of the trace {$\mathbf{P}_2$}--{$P_1$}
              {T}aylor-{H}ood elements for surface {PDE}s},
   JOURNAL = {Math. Comp.},
  FJOURNAL = {Mathematics of Computation},
    VOLUME = {90},
      YEAR = {2021},
    NUMBER = {330},
     PAGES = {1527--1555},
      ISSN = {0025-5718,1088-6842},
   MRCLASS = {65N30 (58J32 65N12 65N15 76D07)},
  MRNUMBER = {4273108},
       DOI = {10.1090/mcom/3551},
       URL = {https://doi.org/10.1090/mcom/3551},
}

@article {TorresSanchezMillanArroyo2019,
    AUTHOR = {Torres-S\'anchez, Alejandro and Mill\'an, Daniel and Arroyo,
              Marino},
     TITLE = {Modelling fluid deformable surfaces with an emphasis on
              biological interfaces},
   JOURNAL = {J. Fluid Mech.},
  FJOURNAL = {Journal of Fluid Mechanics},
    VOLUME = {872},
      YEAR = {2019},
     PAGES = {218--271},
      ISSN = {0022-1120,1469-7645},
   MRCLASS = {76Z05 (92C35)},
  MRNUMBER = {3961595},
       DOI = {10.1017/jfm.2019.341},
       URL = {https://doi.org/10.1017/jfm.2019.341},
}

@article {ArroyoDeSimone2009,
    AUTHOR = {Arroyo, Marino and DeSimone, Antonio},
     TITLE = {Relaxation dynamics of fluid membranes},
   JOURNAL = {Phys. Rev. E (3)},
  FJOURNAL = {Physical Review E. Statistical, Nonlinear, and Soft Matter
              Physics},
    VOLUME = {79},
      YEAR = {2009},
    NUMBER = {3},
     PAGES = {031915, 17},
      ISSN = {1539-3755,1550-2376},
   MRCLASS = {76Z05 (92C10)},
  MRNUMBER = {2497175},
       DOI = {10.1103/PhysRevE.79.031915},
       URL = {https://doi.org/10.1103/PhysRevE.79.031915},
}

@article{PorrmannVoigt2024,
  author  = {Porrmann, Maik and Voigt, Axel},
  title   = {Shape evolution of fluid deformable surfaces under
             active geometric forces},
  journal = {Physics of Fluids},
  volume  = {36},
  number  = {10},
  pages   = {102120},
  year    = {2024},
  doi     = {10.1063/5.0233955}
}

@article{demlow2026taylor,
	title={A Taylor--Hood Finite Element Method for the Surface Stokes Problem Without Penalization},
	author={Demlow, Alan and Neilan, Michael},
	journal={SIAM Journal on Numerical Analysis},
	volume={64},
	number={2},
	pages={565--600},
	year={2026},
	publisher={SIAM}
}

@article{kone2026divergence,
	title={A Divergence-Free Scott-Vogelius Finite Element Method for the Surface Stokes Problem},
	author={Kone, Yerim and Neilan, Michael and Poling, David},
	journal={arXiv preprint arXiv:2606.07840},
	year={2026}
}

@article{fries2018higher,
	title={Higher-order surface FEM for incompressible Navier-Stokes flows on manifolds},
	author={Fries, Thomas-Peter},
	journal={International journal for numerical methods in fluids},
	volume={88},
	number={2},
	pages={55--78},
	year={2018},
	publisher={Wiley Online Library}
}

@article{reusken2025analysis,
	title={Analysis of the Taylor-Hood surface finite element method for the surface Stokes equation},
	author={Reusken, Arnold},
	journal={Mathematics of Computation},
	volume={94},
	number={354},
	pages={1701--1719},
	year={2025}
}

@article{hardering2025parametric,
	title={Parametric finite-element discretization of the surface Stokes equations: inf-sup stability and discretization error analysis},
	author={Hardering, Hanne and Praetorius, Simon},
	journal={IMA Journal of Numerical Analysis},
	volume={45},
	number={5},
	pages={2948--2987},
	year={2025},
	publisher={Oxford University Press}
}

@article{TangWiseLowengrubGuo2026,
	author  = {Tang, Shuqi and Wise, Steven M. and Lowengrub, John S. and Guo, Zhenlin},
	title   = {Membrane heterogeneity-driven dynamics of multicomponent vesicles in shear flow},
	journal = {Journal of Fluid Mechanics},
	year    = {2026},
	volume  = {1030},
	pages   = {A47},
	doi     = {10.1017/jfm.2026.11276}
}

@article {JankuhnOlshanskiiReusken2018,
    AUTHOR = {Jankuhn, Thomas and Olshanskii, Maxim A. and Reusken, Arnold},
     TITLE = {Incompressible fluid problems on embedded surfaces: modeling
              and variational formulations},
   JOURNAL = {Interfaces Free Bound.},
  FJOURNAL = {Interfaces and Free Boundaries. Mathematical Analysis,
              Computation and Applications},
    VOLUME = {20},
      YEAR = {2018},
    NUMBER = {3},
     PAGES = {353--377},
      ISSN = {1463-9963,1463-9971},
   MRCLASS = {35Q30 (37E35 76A20)},
  MRNUMBER = {3875687},
MRREVIEWER = {S\'ergio\ da Silva Rodrigues},
       DOI = {10.4171/IFB/405},
       URL = {https://doi.org/10.4171/IFB/405},
}

@misc{BenavidesNochettoShakipov2025-vector,
    title={$\mathrm{L}^p$-based {S}obolev theory on closed manifolds of minimal regularity: {V}ector-valued problems}, 
    author={Gonzalo A. Benavides and Ricardo H. Nochetto and Mansur Shakipov},
    howpublished={arXiv:2508.11109}, 
    primaryclass = {math.AP},
    year={2025},
}

@article {BenavidesNochettoShakipov2026-scalar,
    AUTHOR = {Benavides, Gonzalo A. and Nochetto, Ricardo H. and Shakipov,
              Mansur},
     TITLE = {{${\rm L}^p$}-based {S}obolev theory on closed manifolds of
              minimal regularity: scalar elliptic equations},
   JOURNAL = {J. Math. Anal. Appl.},
  FJOURNAL = {Journal of Mathematical Analysis and Applications},
    VOLUME = {562},
      YEAR = {2026},
    NUMBER = {1},
     PAGES = {Paper No. 130666, 29},
      ISSN = {0022-247X,1096-0813},
   MRCLASS = {58J05 (35J05 46E35)},
  MRNUMBER = {5060662},
       DOI = {10.1016/j.jmaa.2026.130666},
       URL = {https://doi.org/10.1016/j.jmaa.2026.130666},
}

@Article{BouckNochettoShakipovYushutin2026,
    author={Bouck, Lucas
    and Nochetto, Ricardo H.
    and Shakipov, Mansur
    and Yushutin, Vladimir},
    title={{I}nf-{S}up {S}tability of {P}arabolic {T}race{FEM}},
    journal={Foundations of Computational Mathematics},
    year={2026},
    month={6},
    day={25},
    issn={1615-3383},
    doi={10.1007/s10208-026-09757-7},
    url={https://doi.org/10.1007/s10208-026-09757-7}
}

@article {GrandeLehrenfeldReusken2018,
    AUTHOR = {Grande, J\"org and Lehrenfeld, Christoph and Reusken, Arnold},
     TITLE = {Analysis of a high-order trace finite element method for
              {PDE}s on level set surfaces},
   JOURNAL = {SIAM J. Numer. Anal.},
  FJOURNAL = {SIAM Journal on Numerical Analysis},
    VOLUME = {56},
      YEAR = {2018},
    NUMBER = {1},
     PAGES = {228--255},
      ISSN = {0036-1429,1095-7170},
   MRCLASS = {65N30 (58J32 65N15 65N22)},
  MRNUMBER = {3747467},
MRREVIEWER = {Dietrich\ Braess},
       DOI = {10.1137/16M1102203},
       URL = {https://doi.org/10.1137/16M1102203},
}

@incollection {BonitoDemlowNochetto2020,
    AUTHOR = {Bonito, Andrea and Demlow, Alan and Nochetto, Ricardo H.},
     TITLE = {Finite element methods for the {L}aplace-{B}eltrami operator},
 BOOKTITLE = {Geometric partial differential equations. {P}art {I}},
    SERIES = {Handb. Numer. Anal.},
    VOLUME = {21},
     PAGES = {1--103},
 PUBLISHER = {Elsevier/North-Holland, Amsterdam},
      YEAR = {2020},
      ISBN = {978-0-444-64003-1},
   MRCLASS = {65N30 (35Axx 58J05)},
  MRNUMBER = {4378426},
       DOI = {10.1016/bs.hna.2019.06.002},
       URL = {https://doi.org/10.1016/bs.hna.2019.06.002},
}

@article {Reusken2015,
    AUTHOR = {Reusken, Arnold},
     TITLE = {Analysis of trace finite element methods for surface partial
              differential equations},
   JOURNAL = {IMA J. Numer. Anal.},
  FJOURNAL = {IMA Journal of Numerical Analysis},
    VOLUME = {35},
      YEAR = {2015},
    NUMBER = {4},
     PAGES = {1568--1590},
      ISSN = {0272-4979,1464-3642},
   MRCLASS = {65N30 (65N15)},
  MRNUMBER = {3407236},
MRREVIEWER = {Maxim\ Evgenievich\ Frolov},
       DOI = {10.1093/imanum/dru047},
       URL = {https://doi.org/10.1093/imanum/dru047},
}

@book {GilbargTrudinger2001,
    AUTHOR = {Gilbarg, David and Trudinger, Neil S.},
     TITLE = {Elliptic {P}artial {D}ifferential {E}quations of {S}econd {O}rder},
    SERIES = {Classics in Mathematics},
      NOTE = {Reprint of the 1998 edition},
 PUBLISHER = {Springer-Verlag, Berlin},
      YEAR = {2001},
     PAGES = {xiv+517},
      ISBN = {3-540-41160-7},
   MRCLASS = {35-02 (35Jxx)},
  MRNUMBER = {1814364},
}

@article{JankuhnOlshanskiiReuskenZhiliakov2021,
    AUTHOR  = {Jankuhn, Thomas and Olshanskii, Maxim A. and Reusken,
             Arnold and Zhiliakov, Alexander},
    TITLE   = {Error analysis of higher order trace finite element
             methods for the surface {S}tokes equation},
    JOURNAL = {J. Numer. Math.},
    VOLUME  = {29},
    YEAR    = {2021},
    NUMBER  = {3},
    PAGES   = {245--267},
    DOI     = {10.1515/jnma-2020-0017},
}

@article {PalzhanovZhiliakovQuainiOlshanskii2021,
    AUTHOR = {Palzhanov, Yerbol and Zhiliakov, Alexander and Quaini,
              Annalisa and Olshanskii, Maxim},
     TITLE = {A decoupled, stable, and linear {FEM} for a phase-field model
              of variable density two-phase incompressible surface flow},
   JOURNAL = {Comput. Methods Appl. Mech. Engrg.},
  FJOURNAL = {Computer Methods in Applied Mechanics and Engineering},
    VOLUME = {387},
      YEAR = {2021},
     PAGES = {Paper No. 114167, 21},
      ISSN = {0045-7825,1879-2138},
   MRCLASS = {76T06 (76D05 76M10)},
  MRNUMBER = {4321284},
MRREVIEWER = {Mustafa\ Aggul},
       DOI = {10.1016/j.cma.2021.114167},
       URL = {https://doi.org/10.1016/j.cma.2021.114167},
}

@article{WangPalzhanovQuainiOlshanskiiMajd2022lipid,
 	title={Lipid domain coarsening and fluidity in multicomponent lipid vesicles: A continuum based model and its experimental validation},
 	author={Wang, Y and Palzhanov, Y and Quaini, A and Olshanskii, M and Majd, S},
 	journal={Biochimica et Biophysica Acta (BBA)-Biomembranes},
 	volume={1864},
 	number={7},
 	pages={183898},
 	year={2022},
 	publisher={Elsevier}
}

@article {NitschkeVoigtWensch2012,
    AUTHOR = {Nitschke, I. and Voigt, A. and Wensch, J.},
     TITLE = {A finite element approach to incompressible two-phase flow on
              manifolds},
   JOURNAL = {J. Fluid Mech.},
  FJOURNAL = {Journal of Fluid Mechanics},
    VOLUME = {708},
      YEAR = {2012},
     PAGES = {418--438},
      ISSN = {0022-1120,1469-7645},
   MRCLASS = {76D05 (76M10 76T99)},
  MRNUMBER = {2975450},
       DOI = {10.1017/jfm.2012.317},
       URL = {https://doi.org/10.1017/jfm.2012.317},
}

@article {BonitoDemlowLicht2020,
    AUTHOR = {Bonito, Andrea and Demlow, Alan and Licht, Martin},
     TITLE = {A divergence-conforming finite element method for the surface
              {S}tokes equation},
   JOURNAL = {SIAM J. Numer. Anal.},
  FJOURNAL = {SIAM Journal on Numerical Analysis},
    VOLUME = {58},
      YEAR = {2020},
    NUMBER = {5},
     PAGES = {2764--2798},
      ISSN = {0036-1429,1095-7170},
   MRCLASS = {65N30 (65N12 65N15 65N25 76D07)},
  MRNUMBER = {4155235},
MRREVIEWER = {Barbara\ Verf\"urth},
       DOI = {10.1137/19M1284592},
       URL = {https://doi.org/10.1137/19M1284592},
}

@article {OlshanskiiReuskenGrande2009,
    AUTHOR = {Olshanskii, Maxim A. and Reusken, Arnold and Grande, J\"org},
     TITLE = {A finite element method for elliptic equations on surfaces},
   JOURNAL = {SIAM J. Numer. Anal.},
  FJOURNAL = {SIAM Journal on Numerical Analysis},
    VOLUME = {47},
      YEAR = {2009},
    NUMBER = {5},
     PAGES = {3339--3358},
      ISSN = {0036-1429,1095-7170},
   MRCLASS = {65N30 (65N15 76D45)},
  MRNUMBER = {2551197},
MRREVIEWER = {Riccardo\ Sacco},
       DOI = {10.1137/080717602},
       URL = {https://doi.org/10.1137/080717602},
}

@article {BurmanHansboLarsonMassing2018,
    AUTHOR = {Burman, Erik and Hansbo, Peter and Larson, Mats G. and
              Massing, Andr\'e},
     TITLE = {Cut finite element methods for partial differential equations
              on embedded manifolds of arbitrary codimensions},
   JOURNAL = {ESAIM Math. Model. Numer. Anal.},
  FJOURNAL = {ESAIM. Mathematical Modelling and Numerical Analysis},
    VOLUME = {52},
      YEAR = {2018},
    NUMBER = {6},
     PAGES = {2247--2282},
      ISSN = {2822-7840,2804-7214},
   MRCLASS = {65N30 (58J05 65N15 65N85)},
  MRNUMBER = {3905189},
MRREVIEWER = {Leonardo\ E.\ Figueroa},
       DOI = {10.1051/m2an/2018038},
       URL = {https://doi.org/10.1051/m2an/2018038},
}

@incollection {OlshanskiiReusken2017,
    AUTHOR = {Olshanskii, Maxim A. and Reusken, Arnold},
     TITLE = {Trace finite element methods for {PDE}s on surfaces},
 BOOKTITLE = {Geometrically unfitted finite element methods and
              applications},
    SERIES = {Lect. Notes Comput. Sci. Eng.},
    VOLUME = {121},
     PAGES = {211--258},
 PUBLISHER = {Springer, Cham},
      YEAR = {2017},
      ISBN = {978-3-319-71431-8},
   MRCLASS = {65N30 (35R01)},
  MRNUMBER = {3806652},
       DOI = {10.1007/978-3-319-71431-8\_7},
       URL = {https://doi.org/10.1007/978-3-319-71431-8_7},
}

@misc{Schoberl2009,
	title = {{NGSolve} {Finite} {Element} {Library}},
	url = {http://sourceforge.net/projects/ngsolve},
	urldate = {2022-09-25},
	author = {Schöberl, Joachim},
	month = jan,
	year = {2009},
}

@article{Lehrenfeld2021,
	title = {ngsxfem: {Add}-on to {NGSolve} for geometrically unfitted finite element discretizations},
	volume = {6},
	issn = {2475-9066},
	shorttitle = {ngsxfem},
	url = {https://joss.theoj.org/papers/10.21105/joss.03237},
	doi = {10.21105/joss.03237},
	number = {64},
	urldate = {2022-09-26},
	journal = {Journal of Open Source Software},
	author = {Lehrenfeld, Christoph and Heimann, Fabian and Preuß, Janosch and von Wahl, Henry},
	month = aug,
	year = {2021},
	pages = {3237},
}

@article{BachiniKrauseVoigt2023,
    author = {Bachini, Elena and Krause, Veit and Voigt, Axel},
    title = {The interplay of geometry and coarsening in multicomponent lipid vesicles under the influence of hydrodynamics},
    journal = {Physics of Fluids},
    volume = {35},
    number = {4},
    pages = {042102},
    year = {2023},
    doi = {10.1063/5.0145884}
}

@article {BachiniKrauseNitschkeVoigt2023,
    AUTHOR = {Bachini, Elena and Krause, Veit and Nitschke, Ingo and Voigt,
              Axel},
     TITLE = {Derivation and simulation of a two-phase fluid deformable
              surface model},
   JOURNAL = {J. Fluid Mech.},
  FJOURNAL = {Journal of Fluid Mechanics},
    VOLUME = {977},
      YEAR = {2023},
     PAGES = {Paper No. A41, 37},
      ISSN = {0022-1120,1469-7645},
   MRCLASS = {76A20},
  MRNUMBER = {4681179},
       DOI = {10.1017/jfm.2023.943},
       URL = {https://doi.org/10.1017/jfm.2023.943},
}

@article {OlshanskiiQuainiReuskenYushutin2018,
    AUTHOR = {Olshanskii, Maxim A. and Quaini, Annalisa and Reusken, Arnold
              and Yushutin, Vladimir},
     TITLE = {A finite element method for the surface {S}tokes problem},
   JOURNAL = {SIAM J. Sci. Comput.},
  FJOURNAL = {SIAM Journal on Scientific Computing},
    VOLUME = {40},
      YEAR = {2018},
    NUMBER = {4},
     PAGES = {A2492--A2518},
      ISSN = {1064-8275,1095-7197},
   MRCLASS = {65N30 (65N12 76D07)},
  MRNUMBER = {3841618},
MRREVIEWER = {J\'anos\ Kar\'atson},
       DOI = {10.1137/18M1166183},
       URL = {https://doi.org/10.1137/18M1166183},
}

@article {HeywoodRannacher1982,
    AUTHOR = {Heywood, John G. and Rannacher, Rolf},
     TITLE = {Finite element approximation of the nonstationary
              {N}avier-{S}tokes problem. {I}. {R}egularity of solutions and
              second-order error estimates for spatial discretization},
   JOURNAL = {SIAM J. Numer. Anal.},
  FJOURNAL = {SIAM Journal on Numerical Analysis},
    VOLUME = {19},
      YEAR = {1982},
    NUMBER = {2},
     PAGES = {275--311},
      ISSN = {0036-1429},
   MRCLASS = {65M60 (35Q10 76D05)},
  MRNUMBER = {650052},
       DOI = {10.1137/0719018},
       URL = {https://doi.org/10.1137/0719018},
}

@article {MakridakisNochetto2003,
    AUTHOR = {Makridakis, Charalambos and Nochetto, Ricardo H.},
     TITLE = {Elliptic reconstruction and a posteriori error estimates for
              parabolic problems},
   JOURNAL = {SIAM J. Numer. Anal.},
  FJOURNAL = {SIAM Journal on Numerical Analysis},
    VOLUME = {41},
      YEAR = {2003},
    NUMBER = {4},
     PAGES = {1585--1594},
      ISSN = {0036-1429,1095-7170},
   MRCLASS = {65M15 (35K20)},
  MRNUMBER = {2034895},
MRREVIEWER = {Erik\ Burman},
       DOI = {10.1137/S0036142902406314},
       URL = {https://doi.org/10.1137/S0036142902406314},
}

@misc{OlshanskiiShakipov2026,
    title={Surface {S}tokes-{C}ahn-{H}illiard {S}ystem: {A}nalysis and {S}tructure-{P}reserving {D}iscretization}, 
    author={Maxim Olshanskii and Mansur Shakipov},
    year={2026},
    primaryClass={math.AP},
    howpublished={arXiv:2610.06552},
}

@article {FengProhl2003,
    AUTHOR = {Feng, Xiaobing and Prohl, Andreas},
     TITLE = {Numerical analysis of the {A}llen-{C}ahn equation and
              approximation for mean curvature flows},
   JOURNAL = {Numer. Math.},
  FJOURNAL = {Numerische Mathematik},
    VOLUME = {94},
      YEAR = {2003},
    NUMBER = {1},
     PAGES = {33--65},
      ISSN = {0029-599X,0945-3245},
   MRCLASS = {65M60 (35B25 35K20 35K55 53C44 65M12 65M15 74N99)},
  MRNUMBER = {1971212},
MRREVIEWER = {Stefka\ N.\ Dimova},
       DOI = {10.1007/s00211-002-0413-1},
       URL = {https://doi.org/10.1007/s00211-002-0413-1},
}

@article {FengProhl2004,
    AUTHOR = {Feng, Xiaobing and Prohl, Andreas},
     TITLE = {Error analysis of a mixed finite element method for the
              {C}ahn-{H}illiard equation},
   JOURNAL = {Numer. Math.},
  FJOURNAL = {Numerische Mathematik},
    VOLUME = {99},
      YEAR = {2004},
    NUMBER = {1},
     PAGES = {47--84},
      ISSN = {0029-599X,0945-3245},
   MRCLASS = {65M15 (35B25 35K35 35K55 65M60)},
  MRNUMBER = {2101784},
MRREVIEWER = {Lucia\ Gastaldi},
       DOI = {10.1007/s00211-004-0546-5},
       URL = {https://doi.org/10.1007/s00211-004-0546-5},
}

@article {DziukElliott2013,
    AUTHOR = {Dziuk, Gerhard and Elliott, Charles M.},
     TITLE = {Finite element methods for surface {PDE}s},
   JOURNAL = {Acta Numer.},
  FJOURNAL = {Acta Numerica},
    VOLUME = {22},
      YEAR = {2013},
     PAGES = {289--396},
      ISSN = {0962-4929,1474-0508},
   MRCLASS = {65N30 (58-04)},
  MRNUMBER = {3038698},
MRREVIEWER = {Massimiliano\ Ferronato},
       DOI = {10.1017/S0962492913000056},
       URL = {https://doi.org/10.1017/S0962492913000056},
}

@article {DiegelFengWise2015,
    AUTHOR = {Diegel, Amanda E. and Feng, Xiaobing H. and Wise, Steven M.},
     TITLE = {Analysis of a mixed finite element method for a
              {C}ahn-{H}illiard-{D}arcy-{S}tokes system},
   JOURNAL = {SIAM J. Numer. Anal.},
  FJOURNAL = {SIAM Journal on Numerical Analysis},
    VOLUME = {53},
      YEAR = {2015},
    NUMBER = {1},
     PAGES = {127--152},
      ISSN = {0036-1429,1095-7170},
   MRCLASS = {65M60 (35M33 35Q35 65M12 65M15 76D07)},
  MRNUMBER = {3296618},
MRREVIEWER = {Iuliana\ Stanculescu},
       DOI = {10.1137/130950628},
       URL = {https://doi.org/10.1137/130950628},
}

@article {LiShen2020,
    AUTHOR = {Li, Xiaoli and Shen, Jie},
     TITLE = {On a {SAV}-{MAC} scheme for the
              {C}ahn-{H}illiard-{N}avier-{S}tokes phase-field model and its
              error analysis for the corresponding
              {C}ahn-{H}illiard-{S}tokes case},
   JOURNAL = {Math. Models Methods Appl. Sci.},
  FJOURNAL = {Mathematical Models and Methods in Applied Sciences},
    VOLUME = {30},
      YEAR = {2020},
    NUMBER = {12},
     PAGES = {2263--2297},
      ISSN = {0218-2025,1793-6314},
   MRCLASS = {65M06 (35G25 65M12 65M15 76D07)},
  MRNUMBER = {4179190},
MRREVIEWER = {Reza\ Mohammadi},
       DOI = {10.1142/S0218202520500438},
       URL = {https://doi.org/10.1142/S0218202520500438},
}

@article {ScottZhang1990,
    AUTHOR = {Scott, L. Ridgway and Zhang, Shangyou},
     TITLE = {Finite element interpolation of nonsmooth functions satisfying
              boundary conditions},
   JOURNAL = {Math. Comp.},
  FJOURNAL = {Mathematics of Computation},
    VOLUME = {54},
      YEAR = {1990},
    NUMBER = {190},
     PAGES = {483--493},
      ISSN = {0025-5718,1088-6842},
   MRCLASS = {65D05 (65N30)},
  MRNUMBER = {1011446},
MRREVIEWER = {Qian\ Li},
       DOI = {10.2307/2008497},
       URL = {https://doi.org/10.2307/2008497},
}

\end{document}